%% file: main.tex
\documentclass[webpdf,imanum]{ima-authoring-template}

\usepackage{tikz}
\usetikzlibrary{arrows.meta, decorations.pathmorphing, decorations.markings}
\usepackage{mathrsfs}

\usepackage{pgfplots}
\pgfplotsset{compat=newest}
\usetikzlibrary{plotmarks}
\usepgfplotslibrary{patchplots}
\usepackage{grffile}

\numberwithin{equation}{section}

\theoremstyle{plain}
\newtheorem{theorem}{Theorem}[section]
\newtheorem{lemma}[theorem]{Lemma}
\newtheorem{proposition}[theorem]{Proposition}

\theoremstyle{definition}
\newtheorem{definition}[theorem]{Definition}

\theoremstyle{remark}
\newtheorem{remark}[theorem]{Remark}

\newcommand{\Chat}{\widehat{\mathbb{C}}}

\begin{document}
\DOI{DOI HERE}
\copyrightyear{2025}
\vol{00}
\pubyear{}
\access{Advance Access Publication Date: Day Month Year}
\appnotes{Paper}
\copyrightstatement{Published by Oxford University Press on behalf of the Institute of Mathematics and its Applications. All rights reserved.}
\firstpage{1}

\graphicspath{{figs/}}

%\subtitle{Subject Section}

\title[Low-rank approximation of analytic kernels]{Low-rank approximation of analytic kernels}
\author{Marcus Webb\address{\orgdiv{Department of Mathematics}, \orgname{The University of Manchester}%, \orgaddress{\street{Oxford Road, Manchester} \postcode{M13 9PL}, \country{UK}}
\\ \href{mailto:marcus.webb@manchester.ac.uk}{marcus.webb@manchester.ac.uk}
} {\bf Dedicated to the memory of Nicholas J. Higham}}
\authormark{Marcus Webb}

\received{Date}{0}{Year}
\revised{Date}{0}{Year}
\accepted{Date}{0}{Year}

\editor{Associate Editor:}

\abstract{Many algorithms in scientific computing and data science take advantage of low-rank approximation of matrices and kernels, and understanding why nearly-low-rank structure occurs is essential for their analysis and further development. This paper provides a framework for bounding the best low-rank approximation error of matrices arising from samples of a kernel that is analytically continuable in one of its variables to an open region of the complex plane. Elegantly, the low-rank approximations used in the proof are computable by rational interpolation using the roots and poles of Zolotarev rational functions, leading to a fast algorithm for their construction.}


\maketitle

\section{Introduction}\label{sec:intro}

From system identification \cite{van2012subspace} and particle simulation \cite{greengard1987fast} to image compression \cite{hamlomo2025systematic} and recommender systems \cite{bell2007lessons}, many matrices and kernels in computational science are \emph{nearly-low-rank} --- that is, \emph{near} to low-rank in some norm --- raising the question of why this structure is so common and how it can be exploited in high performance algorithms. 

A widely accepted explanation is that many matrices arise as samples of a smooth kernel, and the approximability of smooth functions by polynomials can be adapted to build low-rank approximations \cite{udell2019big, budzinskiy2025big}. For instance, if $A$ is a matrix with elements $A_{i,j} = K(x_i,y_j)$ for a continuous kernel $K$ and sample points $x_i, y_j \in [-1,1]$, such that for each $x\in [-1,1]$ the function $K(x,\cdot)$ is analytically continuable in the complex plane to a Bernstein ellipse with parameter $R > 1$, then Chebyshev polynomial interpolation yields rank-$n$ matrices $A_n$ satisfying $\|A - A_n\| \leq C R^{-n}$, where $C$ depends on the growth of $K(x,\cdot)$ in the complex plane \cite{little1984eigenvalues, townsend2015continuous}. More recently, \emph{displacement structure} was identified by Beckermann and Townsend as another means of conferring nearly-low-rank structure \cite{beckermann2019bounds, clouatre2025lifting}. Specifically, if $A$ satisfies the Sylvester equation $XA - AY = L$, where $X$ and $Y$ are normal matrices with well separated spectra and $L$ is of low rank, then low-rank approximations to $A$ can be constructed by rational approximation.% functions of $X$ and $Y$. 

This paper develops a new framework for bounding the low-rank approximation error of matrices obtained by sampling a kernel that is analytically continuable in one of its variables from a compact subset of the complex plane to a large open set. The pivotal insight is that such kernels have a latent displacement structure through the Cauchy integral formula. Perhaps surprisingly, these bounds are attainable in practice by rational interpolation with prescribed poles.

To illustrate the superiority of this theory over the current state of the art, consider the matrix,
\begin{equation}\label{eqn:introA}
    A_{i,j} = \frac{\Gamma\left(i+j+\frac{1}{2}\right)}{\Gamma\left(i+j+1\right)}, \qquad i,j=0,\ldots,N,
\end{equation}
where $\Gamma$ is the gamma function \cite[§5]{NIST:DLMF}. This is a positive-definite Hankel matrix that appears within the change-of-basis matrix from Chebyshev coefficients to Legendre coefficients, and a certain FFT-based algorithm for applying the change-of-basis relies on approximation of $A$ by a low-rank matrix \cite{townsend2018fast}. 

\begin{figure}[!t]
\centering
\footnotesize
\input{figs/figure1.tex}
\caption{Comparison of the best low-rank approximation error with various bounds and interpolation errors for the matrix in Equation \eqref{eqn:introA} with $N = 100$. The yellow circles come from an analytically derived Zolotarev rational interpolant based on $[0,N]$ and $[-\infty, -\tfrac12]$ (where $K(x,\cdot)$ is singular for $x \in [0,N]$), whereas the purple triangles come from a numerically derived Zolotarev rational interpolant based on $\{0,1,\ldots, N\}$ and $\{-\tfrac12, -\tfrac32, -\tfrac52, \ldots, -\infty\}$ (where $K(i,\cdot)$ is singular for $i \in \{0,1,\ldots,N\}$).}\label{fig:Hankelintro}
\end{figure}
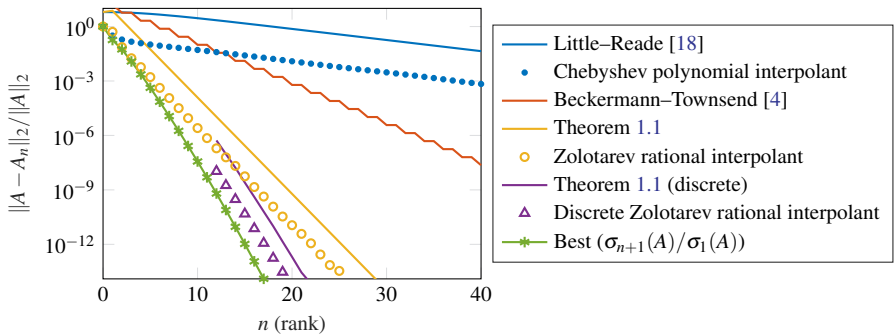

In Figure~\ref{fig:Hankelintro}, we set $N = 100$ and compare the bounds based on polynomial interpolation of $K(x,y) = \Gamma(x+y+1/2)/\Gamma(x+y+1)$ at mapped Chebyshev points on $[0,N]$ (following Little--Reade \cite{little1984eigenvalues}), the bounds of Beckermann--Townsend based on displacement structure implicit in positive-definite Hankel matrices \cite[Cor.~5.5]{beckermann2019bounds} (which ignore any analyticity of the kernel), the bounds and rational interpolants presented in this work, and the singular values of the matrix (the best possible). Further details can be found in Section \ref{sec:BetaMatrices}. Although the existing error bounds from the literature decay exponentially, they are extremely pessimistic compared to Theorem \ref{thm:main}.

The setup of the present work is as follows. Let $E$ be a compact subset of $\mathbb{C}$, let $F$ be a closed subset of $\Chat = \mathbb{C} \cup \{\infty\}$ such that $E \cap F = \emptyset$, and let $D$ be a compact metric space. We assume $D$, $E$, and $F$ are endowed with Radon measures $\lambda$, $\mu$, and $\nu$ respectively. 
Functions on $D$, $E$, and $F$ are denoted by $f$, $g$, and $h$ respectively, and variables in $D$, $E$, and $F$ are denoted by $x$, $y$, and $z$ respectively. 

Any kernel $K \in C(D \times E)$ naturally generates a Hilbert-Schmidt operator $\mathcal{K} : L^2_\mu(E) \to L^2_\lambda(D)$ by
\begin{equation}
    \mathcal{K}[g](x) = \int_E K(x,y) \, g(y) \,\mathrm{d} \mu(y).
\end{equation}
While abstract, this formulation is carefully considered and unifies continuous kernels and discrete matrices and tensors in a single framework. For continuous kernels we take $D$ and $E$ to be continua with $\lambda$ and $\mu$ their respective Lebesgue measures, whereas for a finite dimensional matrix we take
\begin{equation}
    \lambda(x) = \sum_{i=1}^N \delta(x-x_i), \qquad \mu(y) = \sum_{j=1}^M \delta(y-y_j),
\end{equation}
where $\delta$ is the Dirac delta measure, for points $x_1,\ldots,x_N \in D$ and $y_1,\ldots,y_M \in E$, so that $\mathcal{K}$ is represented by the matrix $A \in \mathbb{C}^{N \times M }$, $A_{i,j} = K(x_i, y_j)$, with respect to the Kronecker delta bases (that take the value 1 at a single point and 0 at all others). If $D$ is multidimensional then $A$ is a tensor.

We say that a kernel $K$ has rank $n$ if the range of $\mathcal{K}$ is $n$-dimensional. For a metric space $G$, the space $C(G)$ is that of continuous functions endowed with the topology of boundedness on compact subsets. A function $h$ is analytic at a point $z \in \mathbb{C}$ if there is a Taylor series for $h$ that converges in some disc centered at $z$, and $h$ is analytic at $\infty$ if $h(z^{-1})$ is analytic at $z = 0$. Let $F'$ denote the complement of a subset of $F \subset \Chat$.

\begin{theorem}\label{thm:main}
	Let $K\in C(D \times E)$ be analytically continuable so that $K \in C(D\times F')$ and for each $x \in D$, $K(x,\cdot)$ is analytic in $F'$. Then for $n = 1,2,3,\ldots$, there exists a rank-$n$ kernel $K_n \in C(D \times E)$ such that
	\begin{equation}
		\| \mathcal{K} - \mathcal{K}_n\|_{L_\mu^2(E) \to L_\lambda^2(D)} \leq Z_n\left(L_\mu^2(E),L_\nu^p(F)\right) \, \|\mathcal{K}' \|_{H_{\nu}^p(F)\to L_\lambda^2(D)},
	\end{equation}
where $1 \leq p \leq \infty$, the number $Z_n(L_\mu^2(E),L_\nu^p(F))$ is the Cauchy--Zolotarev number,
    \begin{equation}\label{eqn:CZ}
        Z_n\left(L_\mu^2(E),L_\nu^p(F)\right) = \inf_{\phi \in \mathcal{R}_n} \left\| \frac{\phi(z)^{-1} \phi(y)}{y-z}  \right\|_{L_\mu^2(E) \to L_\nu^p(F)},
    \end{equation}
    where $\mathcal{R}_n$ is the set of all rational functions of type $(n,n)$,
     $\mathcal{K}' : H_\nu^p(F) \to L_\lambda^2(D)$ is the operator
        \begin{equation}\label{eqn:K'}
            \mathcal{K}'[h](x) = \frac{1}{2\pi \mathrm{i}}\int_\Gamma K(x,\xi) \, h(\xi) \, \mathrm{d} \xi,
        \end{equation} 
        defined for any $h$ that is analytic on $F$, which implies that $h$ is analytic in some open set $U \supset F$, and $\Gamma \subset U \setminus F$ is any finite sum of rectifiable Jordan curves with winding numbers satisfying $\mathrm{Ind}_\Gamma(F) - \mathrm{Ind}_\Gamma(E) = 1$. Note that $\mathcal{K}'$ does not depend on the choice of $\Gamma$. The space $H_\nu^p(F) \subseteq L^p_\nu(F)$ is the closure, with respect to the $L^p_\nu(F)$-norm, of the space of functions that are analytic on $F$.
	\end{theorem}
It is the appearance of the numbers $Z_n(L_\mu^2(E),L_\nu^p(F))$ that demonstrates that $\mathcal{K}$ is nearly-low-rank, because these numbers must decay exponentially with $n$ at a rate that depends on how well separated $E$ and $F$ are. Indeed, in Lemma \ref{lem:ZntoZn-1} we show that
\begin{equation}\label{eqn:CZtoZ}
    Z_n\left(L_\mu^2(E),L_\nu^p(F)\right) \leq Z_n(E,F) \, \left\| (y-z)^{-1}  \right\|_{L_\mu^2(E) \to L_\nu^p(F)},
\end{equation}
where
    \begin{equation}\label{eqn:Zolotarev}
	Z_n(E,F) = \min_{\phi \in \mathcal{R}_{n}} \sup_{y \in E, z \in F} |\phi(z)^{-1} \phi(y)|,
\end{equation}
is the $n$th (classical) \emph{Zolotarev number} for $E$ and $F$ \cite{zolotarev1877application, wilber2021computing}. The bound \eqref{eqn:CZtoZ} is not used in the theorem because there are interesting cases in which $\left\| (y-z)^{-1}\right\|_{L_\mu^2(E) \to L_\nu^p(F)} = \infty$ but $Z_n(L_\mu^2(E), L_\nu^p(F))$ decays at the same rate as $Z_n(E,F)$ (see Lemma \ref{lem:Z1}). Section \ref{subsec:logCauchy} contains an example showing that low-rank approximations derived from $Z_n(E,F)$ can be worse than those derived more directly from $Z_n(L_\mu^2(E),L_\nu^p(F))$. The choice of $p$ is therefore delicate. It is often determined by the need to avoid $\|\mathcal{K}' \|_{H_{\nu}^p(F)\to L_\lambda^2(D)}= \infty$, and users of this theorem should start by analysing the action of $\mathcal{K}'$ on analytic functions on $F$, which should elucidate the ideal value of $p$ (usually one of $1$, $2$ or $\infty$).

 It is known that \begin{equation}
    \lim_{n \to \infty} \left( Z_n(E,F) \right)^{1/n} = \exp\left(- \frac{1}{\mathrm{cap}(E,F)}\right) = \rho(E,F)^{-1},
\end{equation} 
where $\mathrm{cap}(E,F)$ is the \emph{logarithmic condenser capacity} of a condenser with plates $E$ and $F$ and $\rho(E,F)$ is known as the \emph{conformal modulus}, both of which are conformal invariants that quantify the separation of subsets of $\Chat$ \cite{gonvcar1969zolotarev}, \cite[Ch.~VII]{nehari1952conformal}. In many cases $Z_n(E,F) = \mathcal{O}(\rho(E,F)^{-n})$ \cite{rubin2022bounding}.

The following is an alternative version of Theorem \ref{thm:main} that is less precise because of the unspecified constant, but requires less machinery to understand. 
\begin{theorem}\label{thm:main2}
    Let $K$, $D$, $E$ and $F$ be as in Theorem \ref{thm:main} and let $1< R < \rho(E,F)$. There exists a constant $C$ such that, for all $n = 1,2,3,\ldots$ there exists a rank-n kernel $K_n \in C(D\times E)$ such that
    \begin{equation*}
        \sup_{x\in D, y \in E}\left|K(x,y) - K_n(x,y)\right| \leq C R^{-n}.
    \end{equation*}
\end{theorem}
The take-home message is that the conformal modulus $\rho(E,F)$ controls the convergence rate of the low-rank approximations and, when $E$ is an interval, is typically far larger than the parameter of the largest Bernstein ellipse, which had previously been the primary quantity used to control this rate in the literature \cite{little1984eigenvalues,townsend2015continuous}.
The conformal modulus is nontrivial to calculate for general sets $E$ and $F$, but if $\Chat \setminus (E\cup F)$ can be conformally mapped to an annulus $\{z \in \mathbb{C} : 1 < |z| < R \}$ (for example, if $E = [-1,1]$ and $F$ is the exterior of a Bernstein ellipse with parameter $R$), then $\rho(E,F) = R$ \cite[Ch.~VII]{nehari1952conformal}. Note that if $E_1 \subseteq E_2$ and $F_1 \subseteq F_2$ then $\rho(E_2,F_2) \leq \rho(E_1,F_1)$, so for a given kernel $K$ one should take $E$ and $F$ as small as possible (to make $E$ and $F$ as well separated as possible).

Equation \eqref{eqn:K'} introduces a new operator $\mathcal{K}'$, which translates the singularities that $K(x,\cdot)$ has in $F$ to an action (such as evaluation, differentiation, or convolution) on functions defined on $F$. Suitable contours $\Gamma$ for general compact sets $F$ and open neighbourhoods $U$ can be constructed as the boundary of a finite union of discs covering $F$, as described in \cite[Sect.~4.2]{grothendieck1953certains}, or using a cartesian grid as in \cite[Sec.~4, Ch.~4]{ahlfors1979complex}. Note that we cannot simply take $\Gamma = \partial F$ because $F$ does not necessarily have regular boundary and even if it does, $K$ is not necessarily bounded on $\partial F$. In Section \ref{sec:Grothendieck} we discuss how $\mathcal{K}'$ can be interpreted in terms of Grothendieck's theory of duality for the function space $H(F; \, L_\lambda^2(D))$ of weakly analytic functions with values in $L_\lambda^2(D)$ \cite{grothendieck1953certains}. In Section \ref{sec:rational}, we show that the low-rank kernel $K_n$ is a rational interpolant of $K(x,\cdot)$, whose interpolation nodes and poles are the roots and poles of the optimal $\phi$ in Equation \eqref{eqn:CZ}. This leads to practical algorithms for the construction of the low-rank approximant by rational interpolation, limited only by the identification of $\phi$. Section \ref{sec:examples} provides several examples of how to apply Theorem \ref{thm:main}, compares the new bounds with the optimal errors, and shows that the computable rational approximations achieve the bounds. Section \ref{sec:conclusion} provides some concluding remarks, and indicates future directions and impact.

\section{Proof of the Main Theorem} 
\begin{figure}[!t]
\centering
\footnotesize
\begin{tikzpicture}[scale=0.85]

  % Styles
  \tikzstyle{blob} = [fill=gray!20, draw=black, thick]
  \tikzstyle{contour} = [thick]
  \tikzstyle{orientarrow} = [postaction={decorate},
    decoration={markings, mark=at position 0.5 with {\arrow{>}}}]
  % --- Draw U (with dotted outline) ---
  \draw[fill=blue!10, draw=black, thick, dotted]
    plot [smooth cycle, tension=1] coordinates {
      (1.3,1.2) (1.3,-0.9) (2.8,-1.7) (4.9,-0.9) (4.8,1.2) (3.1,2.0)
    };
  \node[color=blue!80] at (4.9,0.2) {$U$};

  % Draw blob-like shape for E
  \draw[blob, fill=gray!50]
    plot [smooth cycle, tension=1] coordinates {
      (-3.2,0.3) (-3.4,-0.5) (-2.6,-0.8) (-2.0,-0.2) (-2.3,0.6) (-2.8,0.9)
    }
    node at (-2.8,0.1) {$E$};

  % Draw blob-like shape for F
  \draw[blob, fill=gray!50]
    plot [smooth cycle, tension=1] coordinates {
      (2.7,0.6) (2.4,-0.1) (3.0,-0.7) (3.8,-0.5) (3.6,0.3) (3.2,0.8)
    }
    node at (2.9,0) {$F$};

  % Contour Gamma (around E, clockwise)
  \draw[contour, postaction={decorate}, 
    decoration={markings, mark=at position 0.4 with {\arrow{<}}}]
    plot [smooth cycle, tension=1] coordinates {
      (-3.8,0.7) (-4,-0.6) (-2.5,-1.4) (-1.6,-0.6) (-1.9,0.8) (-3.2,1.6)
    };
  \node at (-4,1.0) {$\tilde\Gamma$};

  % Contour Gamma (around F, counter-clockwise)
  \draw[contour, postaction={decorate},
    decoration={markings, mark=at position 0.4 with {\arrow{>}}}]
    plot [smooth cycle, tension=1] coordinates {
      (2,0.8) (1.8,-0.2) (2.8,-1.2) (4.2,-0.6) (4.0,0.8) (3.0,1.6)
    };
  \node at (1.7,0.9) {$\Gamma$};

  % Point at infinity
  \filldraw[black] (0,1.4) circle (1pt) node[anchor=west, yshift=2pt]{$\infty$};
  \filldraw[black] (-2.35,-0.2) circle (1pt) node[anchor=west, yshift=2pt]{$y$};
  \filldraw[black] (3.4,-0.4) circle (1pt) node[anchor=west, yshift=2pt]{$z$};
  \filldraw[black] (-1.9,0.8) circle (1pt) node[anchor=west, yshift=2pt]{$\tilde\xi$};
  \filldraw[black] (4.0,0.8) circle (1pt) node[anchor=west, yshift=2pt]{$\xi$};

% small circular contour around infinity
\draw[color=gray!300, thick, postaction={decorate},
    decoration={markings, mark=at position 0.25 with {\arrow{>}}}] (0,1.4) circle [radius=0.5];
\node[color=gray!300] at (-0.8,1.4) {$\Gamma_R$}; 

  % Dashed deformation arrow with label
  \draw[->, thick, dashed, bend left=-20, color=gray!300, font=\footnotesize] (-1.45,0) to node[below] {deformation} (1.65,0);
  \node[text=gray!300, right, font=\footnotesize] at (-2.4,2.2) {zero contribution};
\end{tikzpicture}
\caption{A diagram showing the relationship between the sets, contours and points occurring in the proof of Theorem \ref{thm:main}. Here $\mathrm{Ind}_{\tilde{\Gamma}}(E) = -1$, $\mathrm{Ind}_{\tilde{\Gamma}}(F) = 0$, $\mathrm{Ind}_{\Gamma}(E) = 0$, $\mathrm{Ind}_{\Gamma}(F) = 1$. Note that $\infty$ could also be located within $F$, discussed in the text. }\label{fig:maintheorem}
\end{figure}
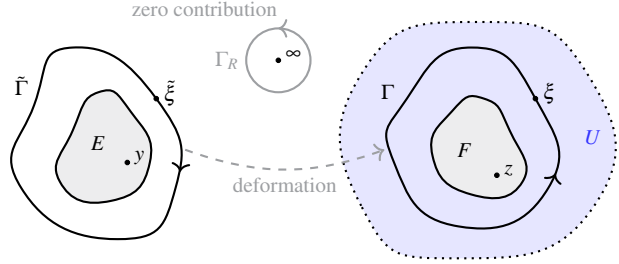

By Cauchy's integral formula, for any $x \in D$, $y\in E$,
\begin{equation}
    K(x,y) = \frac{1}{2 \pi \mathrm{i}}\int_{\tilde\Gamma} \frac{K(x,\tilde\xi)}{y-\tilde\xi} \, \mathrm{d} \tilde\xi,
\end{equation}
where $\tilde\Gamma \subset F'\setminus E$ is any finite sum of rectifiable Jordan curves, with winding numbers satisfying $\mathrm{Ind}_{\tilde\Gamma}(F) - \mathrm{Ind}_{\tilde\Gamma}(E) = 1$ (note that this goes against the usual sign convention). If $\infty \in F$ then for any open set $U \supset F$, we have that $U'$ is a bounded set containing $E$ and we can deform the contour to $\Gamma \subset U\setminus F$ and preserve the winding numbers. If $F$ is a compact set as in Figure \ref{fig:maintheorem}, so that $\infty \in \Chat \setminus (E\cup F)$, then we must take care when deforming over the point $\infty$. The decomposition $K(x,\xi)(y-\xi)^{-1} = a(x,y)\xi^{-1} + b(x,y,\xi)\xi^{-2}$, valid for $|\xi| > R$ for sufficiently large $R$, where $b(x,y,\xi)$ is bounded as $\xi \to \infty$ can be used to show that the contribution of the contour $\Gamma_R = \{\xi: |\xi| = R \}$ surrounding $\infty$ is zero, leaving the only contributing contour to be $\Gamma$ with modified winding numbers (because of the dropped contour around infinity), but still satisfying $\mathrm{Ind}_\Gamma(F) - \mathrm{Ind}_\Gamma(E)= 1$.

Now, for $g \in L_\mu^2(E)$, by Fubini's Theorem,
\begin{equation}\label{eqn:firstdecomp}
\mathcal{K}[g](x) = \frac{1}{2\pi \mathrm{i}}  \int_E \left(\int_\Gamma \frac{K(x,\xi)}{y-\xi} \, \mathrm{d} \xi \right)\, g(y) \, \mathrm{d}\mu(y)  = \frac{1}{2\pi \mathrm{i}}   \int_\Gamma K(x,\xi) \, \left( \int_E\frac{g(y)}{y-\xi}\, \mathrm{d}\mu(y) \right) \, \mathrm{d} \xi.
\end{equation}
Application of Fubini's Theorem is justified because the integrand is absolutely integrable on the compact set $E \times \Gamma$ (as $L_\mu^2(E) \subseteq L_\mu^1(E)$). From this we are led to the key factorisation of $\mathcal{K}$. But first we must define the appropriate spaces for the operators in this factorisation. 
\begin{definition}\label{def:HU}
Let $U$ be an open subset of $\Chat$. The set $H(U)$ consists of all analytic functions $h$ on $U$ such that if $\infty \in U$ then $h(\infty) = 0$. Let $V$ be an arbitrary subset of $\Chat$. The set $H(V)$ consists of equivalence classes of functions (germs), each of which are members of $H(U)$ for some open set $U \supset V$.
\end{definition}
\begin{lemma}\label{lem:factorisation}
Let $K$ be as in Theorem \ref{thm:main}. Then
    \begin{equation}
        \mathcal{K} = \mathcal{K}' \, \mathcal{C}, \qquad \mathcal{K} : L_\mu^2(E) \overset{\mathcal{C}}\to H(F) \overset{\mathcal{K}'}\to C(D),
    \end{equation}
    where $\mathcal{K}'$ is as in Equation \eqref{eqn:K'} and $\mathcal{C}$ is the Cauchy transform
\begin{equation}\label{eqn:Cauchytransform}
    \mathcal{C}[g](z) = \int_E \frac{g(y)}{y - z} \, \mathrm{d} \mu(y).
\end{equation}%
\begin{proof}
    First we prove that $\mathcal{K}'$ and $\mathcal{C}$ map to and from the stated spaces. By differentiating under the integral sign, the Cauchy transform maps to an analytic function on $E'$. Furthermore, $\mu(E) < \infty$ (since $\mu$ is locally finite and $E$ is bounded). Therefore, $\lim_{z \to \infty} \mathcal{C}[g](z) = 0$, so
    $\mathcal{C} : L_\mu^2(E) \to H(F)$.

Let $h \in H(F)$. We need to show that $f = \mathcal{K}'[h] \in C(D)$. Since $F$ is a closed subset of $\Chat$, there exists an open set $U$ containing $F$ such that $h \in H(U)$. Further, for any $x \in D$, $K(x,\cdot) \in C(F')$, so both $K(x,\xi)$ and $h(\xi)$ are bounded for $\xi \in \Gamma$. Therefore, $f(x)$ is defined for all $x \in D$. Now consider a sequence $\{x_k\}_{k=1}^\infty \subset D$ converging to $x \in D$, and associated with that, consider the sequence of functions $ K(x_k,\cdot) \in C(F')$. Since $K \in C(D\times F')$ and by the Heine--Cantor Theorem (continuous functions on compact sets are uniformly continuous), we have that $K(x_k,\cdot) \to K(x,\cdot)$ uniformly on compact subsets of $F'$. Since $\Gamma$ is a compact subset of $F'$, we have $f(x_k) = \int_\Gamma K(x_k, \xi) h(\xi) \,\mathrm{d}\xi \to \int_\Gamma K(x, \xi) h(\xi) \,\mathrm{d}\xi = f(x)$ as $k \to \infty$. Therefore, $f$ is continuous on $D$.

Now let us prove that $\mathcal{K} = \mathcal{K}' \, \mathcal{C}$. By Equation \eqref{eqn:firstdecomp}, we have $\mathcal{K}[g](x) = \frac{1}{2\pi \mathrm{i}}   \int_\Gamma K(x,\xi) \, \mathcal{C}[g](\xi) \, \mathrm{d} \xi$. Note that we have evaluated $\mathcal{C}[g]$ along $\Gamma$, which is outside of $F$. This is justified because $\mathcal{C}[g] \in H(F)$, so it is also in $H(U)$ for some open $U \supset F$ (namely, $U = E'$).
\end{proof}
\end{lemma}
\begin{lemma}[{cf.~\cite{tyrtyshnikov1996mosaic}, \cite[Rm.~4.3]{beckermann2019bounds}}]\label{lemma:C_n}
The kernel $C_n \in C(F \times E)$,
\begin{equation}\label{eqn:Cn}
    C_n(z,y) = \left( 1 - \frac{\phi(y)}{\phi(z)} \right) \, \frac{1}{y-z},
\end{equation}
is of rank at most $n$ for any $\phi \in \mathcal{R}_n$ such that $\phi$ is bounded on $E$ and $1/\phi$ is bounded on $F$.
\begin{proof}
    Write $\phi(\xi) = p(\xi)/q(\xi)$, where $p$ and $q$ are polynomials of degree at most $n$. Then
    \begin{equation}
        C_n(z,y) =  \frac{(y-z)^{-1}(p(z)q(y) - p(y)q(z))}{p(z)q(y)}.
    \end{equation}
    For a fixed $z\in F \setminus \{\infty\}$, we deduce that $C_n(z,\cdot) \in \mathcal{R}_{n-1,n}$, because the denominator is a polynomial of degree at most $n$ and the numerator is a polynomial of degree at most $n$ with a root at $y = z$ that has been factored out to yield a polynomial of degree at most $n-1$. If $\infty \in F$, then since $1/\phi$ is bounded on $F$, we have that $C_n(\infty, y) = 0$ for all $y \in E$, so $C_n(\infty,\cdot) \in \mathcal{R}_{n-1,n}$ trivially. Expansion in partial fractions establishes that the rank is at most $n$.
    % Now, let us write $q_1,\ldots,q_{\ell} \in \mathbb{C}$ be the distinct zeros of $\phi$ and $p_1,\ldots, p_{m} \in \mathbb{C}$ for the zeros and poles of $\phi$, each appearing as many times as their multiplicity. By inspecting \eqref{eqn:Cn} directly, we see that $C_n(z,\cdot)$ interpolates $C(z,\cdot)$ at $q_1,\ldots,q_{\ell}$ (in the Hermite sense if necessary), and has poles at $p_1,\ldots,p_{m}$. Therefore, $C_n(z, \cdot)$ is the unique rational interpolant, and can be written in Lagrange form,
    %\begin{equation}
    %    C_n(z,y) =  \sum_{j=1}^{\ell} C(z,q_j) \phi_j(y), \qquad \phi_j(y) = \left(\prod_{\substack{k=1 \\ k \neq j}}^\ell \frac{ y - q_k}{q_j-q_k} \right) \cdot \left(\prod_{k=1}^m\frac{q_j - p_k}{y-p_k} \right).
    %\end{equation}. 
\end{proof}
\end{lemma}

As long as $\phi$ is bounded on $E$ and $1/\phi$ is bounded on $F$, the associated operator $\mathcal{C}_n$ maps $L_\mu^2(E)$ to $H(F)$ just as in Lemma \ref{lem:factorisation} (except this time $\mathcal{C}_n[g]$ may not be analytic in all of $E'$ because of the poles of $
\phi$). Therefore, we can define $\mathcal{K}_n = \mathcal{K}' \, \mathcal{C}_n$ as a mapping from $L_\mu^2(E)$ to $C(D)$. This operator has a range with dimension at most $n$, because $\mathcal{C}_n$ does, so it is of rank at most $n$. Furthermore,
\begin{equation}\label{eqn:compbound}
    \|\mathcal{K}-\mathcal{K}_n\|_{L_\mu^2(E) \to L_\lambda^2(D)} \leq \|\mathcal{K}'\|_{H_\nu^p(F) \to L_\lambda^2(D)} \, \|\mathcal{C} - \mathcal{C}_n\|_{L_\mu^2(E) \to L_\nu^p(F)},
\end{equation}
for any space $H_\nu^p(F)$ for which these norms are finite. The operator $\mathcal{C}-\mathcal{C}_n$ has kernel
\begin{equation}\label{eqn:CminusC_n}
    C(z,y) - C_n(z,y) = \frac{\phi(z)^{-1}\phi(y)}{y-z}.
\end{equation}
From this it follows that $\inf_{\phi \in \mathcal{R}_n} \|\mathcal{C}-\mathcal{C}_n\|_{L_\mu^2(E) \to L_\nu^p(F)} = Z_n(L_\mu^2(E), L_\nu^p(F))$, from which Theorem \ref{thm:main} follows after substitution into Equation \eqref{eqn:compbound}.

Let us finish this section with a proof of Theorem \ref{thm:main2}. Let $1 < R < \rho(E,F)$. By taking the closure of a finite union of sufficiently small open discs covering $F$, there exists a closed set $G \supset F$ with the following properties: $\Gamma := \partial G$ is a finite union of rectifiable Jordan curves upon which $K(x,\cdot)$ is bounded, $R < \rho(E,G) < \rho(E,F)$, and  $\min_{y\in E, \xi \in \Gamma}|z-\xi| > 0$. By Lemma \ref{lem:factorisation} and Equation \eqref{eqn:CminusC_n},
\begin{equation}
    K(x,y) - K_n(x,y) = \frac{1}{2\pi\mathrm{i}} \int_\Gamma K(x,\xi) \frac{\phi(y)}{\phi(\xi)}\frac{1}{y - \xi} \, \mathrm{d}\xi,
\end{equation}
for any $x \in D$, $y \in E$. We are free to choose $\phi \in \mathcal{R}_n$. Therefore, by bounding the contour integral in terms of the maximum absolute value of the integrand,
\begin{equation}
    |K(x,y) - K_n(x,y)| \leq \frac{\max_{x \in D, \xi \in\Gamma} |K(x,\xi)|}{2\pi\min_{y\in E, \xi \in \Gamma} |y-\xi|}  \, \cdot \, \min_{\phi\in\mathcal{R}_n} \max_{y\in E,\xi\in\Gamma}|\phi(\xi)^{-1}\phi(y)|.
\end{equation}
The first multiplicand is a constant independent of $n$. The second multiplicand is bounded by
\begin{equation}
    \inf_{\phi\in\mathcal{R}_n} \max_{y\in E,\xi\in\Gamma}|\phi(\xi)^{-1}\phi(y)| = \inf_{\phi\in\mathcal{R}_n} \max_{y\in E,\xi\in G}|\phi(\xi)^{-1}\phi(y)| = Z_n(E,G).
\end{equation}
Now, note that $\lim_{n\to \infty} (Z_n(E,G))^{1/n} = \rho(E,G)^{-1} < R^{-1}$. This implies that $Z_n(E,G) = \mathcal{O}(R^{-n})$.

\section{Cauchy--Zolotarev numbers}\label{sec:Cauchy-Zolotarev}

The Cauchy--Zolotarev numbers, a new concept, provide the rapid decay in Theorem \ref{thm:main} and are defined as follows.
\begin{definition}[Cauchy--Zolotarev Numbers]
Let $E$ and $F$ be disjoint, closed sets in $\Chat$, with Radon measures $\mu$ and $\nu$, and let $1 \leq p, q \leq \infty$. Then for $n=0,1,2,\ldots$, the $n$th Cauchy--Zolotarev number is
\begin{equation}
    Z_n\left(L_\mu^p(E),L_\nu^q(F)\right) = \inf_{\phi \in \mathcal{R}_n} \left\|\frac{\phi(z)^{-1} \phi(y)}{y-z} \right\|_{L_\mu^p(E) \to L_\nu^q(F)},
    \end{equation}
    where $\mathcal{R}_n = \mathcal{R}_{n,n}$ is the set of all rational functions of type $(n,n)$.
\end{definition}
These numbers are related to both the classical Zolotarev numbers in Equation \eqref{eqn:Zolotarev} and the Cauchy transform in Equation \eqref{eqn:Cauchytransform}. Indeed, the case $n = 0$ is, $Z_0(L_\mu^p(E), L_\nu^q(F)) = \|\mathcal{C}\|_{L_\mu^p(E)\to L_\nu^q(F)}$, and for $n \geq 1$, the relationship is encapsulated in the following lemma.
\begin{lemma}\label{lem:ZntoZn-1}
Let $n \in \mathbb{Z}_{\geq 0}$ and $ s \in \{0,1,\ldots, n\}$. Then
    \begin{equation}
        Z_n(L_\mu^p(E), L_\nu^q(F)) \leq Z_{n-s}(E,F) \, \cdot \, Z_s(L_\mu^p(E), L_\nu^q(F)).
    \end{equation}
    \begin{proof}
    Consider $\phi(\xi) = \phi_{n-s}(\xi) \phi_s(\xi)$, where $\phi_{n-s} \in \mathcal{R}_{n-s}$, $\phi_s \in \mathcal{R}_s$. Then
    \begin{eqnarray*}
        \left\| \frac{\phi(z)^{-1} \phi(y)}{y-z} \right\|_{L_\mu^p(E) \to L_\nu^q(F)} &=& \sup_{\|g\|_p = 1} \left\| \int_E \frac{\phi(z)^{-1} \phi(y) g(y)}{y-z} \, \mathrm{d}\mu(y) \right\|_{L_\nu^q(F)} \\
        &\leq& \sup_{\|g\|_p = 1} \|\phi_{n-s}^{-1}\|_{L_\nu^\infty(F)}\,\left\| \int_E \frac{ \phi_s(z)^{-1} \phi_s(y) \left(\phi_{n-s}(y) g(y)\right)}{y-z} \, \mathrm{d}\mu(y) \right\|_{L_\nu^q(F)} \\
        &\leq& \sup_{\|g\|_p = 1} \|\phi_{n-s}^{-1}\|_{L_\nu^\infty(F)} \, \left\| \frac{ \phi_s(z)^{-1} \phi_s(y)}{y-z}\right\|_{L_\mu^p(E) \to L_\nu^q(F)} \left\|\phi_{n-s}(y) g(y)\right\| \\
        &\leq& \|\phi_{n-s}^{-1}\|_{L_\nu^\infty(F)}\, \|\phi_{n-s}\|_{L_\mu^\infty(E)}\, \left\| \frac{ \phi_s(z)^{-1} \phi_s(y)}{y-z} \right\|_{L_\mu^p(E) \to L_\nu^q(F)}.
        \end{eqnarray*}
        Now if we take the infimum over all possible $\phi_{n-s}$ and $\phi_s$, we obtain the desired inequality.
    \end{proof}
\end{lemma}
The case $s = 0$ corresponds to the bound
\begin{equation}\label{eqn:ZtoC}
    Z_n(L_\mu^p(E), L_\nu^q(F)) \leq Z_{n}(E,F) \,  \|\mathcal{C}\|_{L_\mu^p(E) \to L_\nu^q(F)}.
\end{equation}
This is sufficient to show that the Cauchy--Zolotarev numbers decay extremely quickly if $\mathcal{C}$ is bounded, because so do the classical Zolotarev numbers. Following the brief discussion in the introduction, let us discuss two important special cases for $E$ and $F$. If $E$ and $F$ are two disjoint disks in the complex plane, then $Z_n(E, F) = \rho(E,F)^{-n}$ \cite{starke1992near}. In particular, for $R > 1$,
    \begin{equation}
	    Z_n(\{|z| \leq 1 \}, \{|z| \geq R \}) = R^{-n}.
	\end{equation}
The Zolotarev numbers are invariant under M\"obius transformations, so the result extends to any two disjoint disks. The rate of decay can be explicitly bounded for the case of two disjoint intervals, too.
% \begin{theorem}[Beckermann--Townsend\cite{beckermann2019bounds}]\label{thm:bt}
% 	Let $E = [a,b]$, $F=[c,d]$, where $-\infty \leq a < b < c < d \leq \infty$. Then for $n = 1,2,3,\ldots$,
% 	\begin{equation}
% 		Z_n(E,F) = 4 \rho^{-n} \left(\prod_{\tau=1}^\infty \frac{1+\rho^{-4\tau n}}{1+\rho^{(2-4\tau)n}} \right)^4, \qquad \rho = \exp\left(\frac{\pi^2}{2\mu(1/\sqrt{\gamma})} \right),
% 	\end{equation}
% 	where
%     \begin{equation}
%         \mu(\lambda) = \frac{\pi}{2}\frac{K(\sqrt{1-\lambda^2})}{K(\lambda)}, \qquad \gamma = \frac{|c-a||d-b|}{|c-b||d-a|},
%     \end{equation}
%     Here $\mu$ is called the Gr\"otzsch ring function and $\gamma$ is the cross ratio. Not to be confused with the kernels in this paper, $K$ is the complete elliptic integral of the first kind.
% \end{theorem}

\begin{theorem}[Beckermann--Townsend \cite{beckermann2019bounds}]\label{thm:BT}
Let $-\infty \leq a < b < c < d \leq \infty$. Then for $n = 1,2,\ldots$,
\begin{equation*}
    Z_n([a,b],[c,d]) \leq 4 \exp\left(  \frac{-n\pi^2}{\log(16 \gamma)}\right), \quad \text{ where } \quad \gamma = \frac{(c-a)(d-b)}{(c-b)(d-a)}. 
\end{equation*}
\end{theorem}

The Cauchy--Zolotarev numbers can be bounded using double integrals, as follows. For compactness of notation, and because we will use these spaces later, we use the Bochner spaces, $L_\mu^p\left(E; \, L_\nu^q(F)\right)$, of Bochner measurable functions $Z$ on $E$ with values in $L_\nu^q(F)$, satisfying finiteness of the norm,
\begin{equation}
    \|Z\|_{L_\mu^p\left(E; \, L_\nu^q(F)\right)} =  \left( \int_E \left( \int_F |Z(z,y)|^q \, \mathrm{d}\nu(z) \right)^{p/q} \,\mathrm{d}\mu(y) \right)^{1/p},
\end{equation}
in the cases $1 \leq p,q < \infty$ and appropriate suprema when $p$ or $q$ is infinite \cite{diestel1977vector}.
\begin{lemma}\label{lem:boundZ}
Let $Z \in C(F \times E)$ define an operator $\mathcal{Z} : L_\mu^{p^*}(E) \to L_\nu^q(F)$ for $p,q \in [1,\infty]$, where $p^*$ is the H\"older dual exponent, satisfying $\tfrac{1}{p} + \tfrac{1}{p^*} = 1$. Then
\begin{equation*}
    \|\mathcal{Z}\|_{L_\mu^{p^*}(E) \to L_\nu^q(F)} \leq \|Z\|_{L_\mu^p\left(E; \, L_\nu^q(F)\right)}.
\end{equation*}
\begin{proof}
    We can work through this directly from the definition. For $p \in [1,\infty]$ and $q \in [1,\infty)$,
        \begin{eqnarray*}
            \|\mathcal{Z}\|_{L_\mu^{p^*}(E) \to L_\nu^q(F)} &=& \sup_{\|g\|_{p^*}=1} \left(\int_F \left| \int_E Z(z,y)\, g(y) \, \mathrm{d}\mu(y) \right|^q \, \mathrm{d} \nu(z)\right)^{1/q} \\
            &\leq& \sup_{\|g\|_{p^*}=1} \int_E \left(\int_F \left|Z(z,y) \right|^q \, \mathrm{d}\nu(z)  \right)^{1/q} \, |g(y)| \, \mathrm{d}\mu(y) \qquad \text{ (Minkowski's inequality)}.
        \end{eqnarray*}
        The result now follows by H\"older's inequality. For $q = \infty$, the proof is similar.
    \end{proof}
\end{lemma}
The Cauchy transform operator is not necessarily bounded for any given $p$ and $q$. The choice $p = 1$, $q =\infty$ always yields bounded $\mathcal{C}$ as long as $E \cap F = \emptyset$, in which case Lemma \ref{lem:boundZ} implies
\begin{equation}
    Z_n(L_\mu^1(E), L_\nu^\infty(F)) = \sup_{y\in E,z \in F} \left| \frac{\phi(z)^{-1}\phi(y)}{z-y} \right| \leq \frac{Z_n(E,F)}{\inf_{y\in E, z \in F} |z-y|} < \infty.
\end{equation}
Some elegant special cases can be derived from Equation \eqref{eqn:ZtoC} and Lemma \ref{lem:boundZ}. For example, with the notation of Theorem \ref{thm:BT}, we can show that $\|\mathcal{C}\|_{L^2(a,b)\to L^2(c,d)} \leq \log^{1/2}(\gamma)$, so
\begin{equation}
     Z_n(L^2(a,b), L^2(c,d)) \leq 4\,\log^{1/2}(\gamma) \, \exp\left(  \frac{-n\pi^2}{\log(16 \gamma)}\right).
\end{equation}
However, in the context of Theorem \ref{thm:main}, we are often constrained by $\mathcal{K}'$ in our choice of $p$ and $q$, and in some of the examples of Section \ref{sec:examples},  $\mathcal{C}$ must be unbounded. This is fatal for the case $n=0$, but does not affect the validity of the main result for $n\geq 1$, as shown in the following Lemma.
\begin{lemma}\label{lem:Z1}
Let $-\infty < b < c <  d< \infty$ and let $\mu$ be a Radon measure on $[c,d]$. Then $Z_0(L_\mu^2(c,d), L^1(-\infty,b)) = \infty$, but for $n \geq 1$,
    \begin{equation*}
        Z_n(L_\mu^2(c,d), L^1(-\infty,b)) \leq  Z_{n-1}([c,d],[-\infty,b]) \, \frac12 \, \mu([c,d])^{1/2} \, \log\left(\frac{d-b}{c-b} \right).
    \end{equation*}
    Note that $(d-b)/(c-b)$ is $\gamma$ in the limiting case $a \to -\infty$.
    \begin{proof}
    For the first assertion, we note that,
\begin{eqnarray}
    Z_0(L_\mu^2(c,d),L^1(-\infty,b)) &=& \sup_{\|g\|_2 = 1} \int_{-\infty}^b \left| \int_c^d \frac{g(y)}{y-z} \, \mathrm{d}\mu(y)\right| \, \mathrm{d}z \\
    &\geq& \mu([c,d])^{-1/2}\int_{-\infty}^b \int_c^d \frac{1}{y-z} \, \mathrm{d}\mu(y) \, \mathrm{d}z,
\end{eqnarray}
by taking $g \equiv \mu([c,d])^{-1/2}$ and noting that $y > z$ in the integrand. Then by Fubini's Theorem,
\begin{equation}
Z_0(L_\mu^2(c,d),L^1(-\infty,b)) \geq \mu([c,d])^{-1/2}\int_c^d \int_{-\infty}^b \frac{1}{y-z} \,\mathrm{d}z \, \mathrm{d}\mu(y) \geq \mu([c,d])^{1/2} \int_{-\infty}^b \frac{1}{d-z}\,\mathrm{d} z = \infty.
\end{equation}
    For the second assertion, choose $\phi(\xi) = \xi-t$, where $t=b + \sqrt{(d-b)(c-b)}$. If $Z(z,y) = \phi(z)^{-1} (y-z)^{-1} \phi(y)$, then by Lemma \ref{lem:boundZ} with $p = 2$, $q = 1$, we have 
        \begin{eqnarray*}
          Z_1(L_\mu^2(c,d), L^1(-\infty,b)) &\leq&
           \left( \int_c^d \left(\int_{-\infty}^b |Z(z,y)| \, \mathrm{d} z \right)^2 \, \mathrm{d}\mu(y) \right)^{1/2}.
        \end{eqnarray*}
         Then, using partial fractions to evaluate the inner integral,
        \begin{eqnarray}
\left( \int_c^d \left(\int_{-\infty}^b |Z(z,y)| \, \mathrm{d} z \right)^2 \, \mathrm{d}\mu(y) \right)^{1/2} &\leq& \left( \int_c^d \left(\int_{-\infty}^b \frac{y-t}{(t-z)(y-z)} \, \mathrm{d} z \right)^2 \, \mathrm{d}\mu(y) \right)^{1/2} \\
&=& \left( \int_c^d \log^2\left( \frac{y-b}{t-b} \right) \, \mathrm{d}\mu(y) \right)^{1/2} \\
&\leq& \mu([c,d])^{1/2} \sup_{y \in [c,d]} \left|\log\left( \frac{y-b}{t-b} \right)\right|.
        \end{eqnarray}
        Now we can substitute the form of $t$ into the logarithm. Since $y \mapsto \log\left(\frac{y-b}{t-b} \right)$ is a monotonic function, its extrema are attained at the endpoints. For this choice of $t$, the endpoint values are $\pm \frac12 \log \left(\frac{d-b}{c-b} \right)$,
        which yields the desired bound on $Z_1$. To bound $Z_n$ for $n > 1$, use Lemma \ref{lem:ZntoZn-1} with $s = 1$.
    \end{proof}
\end{lemma}

%If one finds oneself in a situation in which the Cauchy transform is not a bounded operator, then besides changing $p$, $q$, one could also weight the spaces with a choice of $\mu$ or $\nu$.

\section{The Grothendieck dual operator}\label{sec:Grothendieck}

The operator $\mathcal{K}' : H(F) \to C(D)$ in Theorem \ref{thm:main} appears to be an unusual choice, but we will see here how it can be thought of naturally as a kind of dual operator to $\mathcal{K}$. The type of duality required is due to Grothendieck \cite{grothendieck1953certains}, K\"othe \cite{kothe1953dualitat} and Silva \cite{silva1950funcoes}, and does not appear to be well known among Numerical Analysts, despite the attempt made in \cite{rubel1969functional} to introduce it to our repertoire. %\footnote{The original paper, in Portugese, can be found (as of 2025) at \href{https://purl.pt/2193}{https://purl.pt/2193}}.
We refer to this duality as Grothendieck duality and refer to $\mathcal{K}'$ as the \emph{Grothendieck dual} of $\mathcal{K}$.

\subsection{Examples}

To motivate Grothendieck duality, we begin with two examples. The first is
\begin{equation}
K(x,y) = \frac{m!}{(y-x)^{m+1}}, \qquad x \in D = [1,N], \quad y \in E = [-N,-1]
\end{equation}
where $m$ is an integer. For each $x \in D$, $K(x,\cdot)$ is analytic in $F'$ where $F = [1,N]$. Its Grothendieck dual operator $\mathcal{K}' : H([1,N]) \to C([1,N])$ is
\begin{equation}
    \mathcal{K}'[h](x) = \frac{m!}{2\pi\mathrm{i}}\int_\Gamma \frac{h(\xi)}{(\xi-x)^{m+1}} \, \mathrm{d} \xi = h^{(m)}(x),
\end{equation}
for $h \in H([1,N])$. The second example is
\begin{equation}\label{eqn:logkerneldef}
K(x,y) =\log(y-x), \qquad x \in D = [-N,-1], \quad y \in E=[1,N].
\end{equation}
For each $x \in D$, $K(x,\cdot)$ is analytic in $F'$ where $F = [-\infty, -1]$. The action of the the Grothendieck dual operator $\mathcal{K}' : H([-\infty,-1]) \to C([-N,-1])$ is given below, but proved in Section \ref{sec:proofs}.
\begin{lemma}\label{lem:logkernel}
Let $K$ be as in Equation \eqref{eqn:logkerneldef} and $h \in H([-\infty,-1])$, then for $x \in [-N,-1]$,
    \begin{equation}
        \mathcal{K}'[h](x) = \int_{-\infty}^{x} \frac{a_1}{z-x-1} - h(z) \, \mathrm{d} z,
    \end{equation}
    where $a_1 = \lim_{z \to \infty} z h(z)$.
\end{lemma}

\subsection{Grothendieck duality}
Recall Definition \ref{def:HU}, which defined the spaces $H(F)$ and $H(F')$. Endow $H(F')$ with the topology of uniform convergence on compact subsets of $F'$, and endow $H(F)$ with the inductive limit topology, which means that $h_m \to h$ in $H(F)$ as $m \to \infty$ if there exists an open set $U \supset F$ such that, for all sufficiently large $m$, $h_m \in H(U)$ and $h_m \to h$ in $H(U)$. 

Let $H'(F)$ be the space of continuous linear functionals on $H(F)$. Silva proved that $H'(F)$ is isomorphic to $H(F')$ (as topological vector spaces, if $H'(F)$ is endowed with the strong topology) by the map $\Phi : H(F') \to H'(F)$,
\begin{equation}
    \left(\Phi[g]\right)[h] = \frac{1}{2\pi \mathrm{i}} \int_\Gamma g(\xi)\,h(\xi) \, \mathrm{d}\xi, \qquad g \in H(F'),\quad  h \in H(F),
\end{equation}
where $h$ is analytic in an open set $U$ containing $F$ and $\Gamma \subset U \setminus F$ is any finite sum of rectifiable Jordan curves with winding numbers satisfying $\mathrm{Ind}_\Gamma(F) - \mathrm{Ind}_\Gamma(U')= 1$ \cite{silva1950funcoes}. The inverse operator is $\Phi^{-1}[\varphi](y) = \varphi[(y-z)^{-1}]$, where $\varphi \in H'(F)$ and $y \in F'$. If $F$ is a smooth contour then $\Phi[g]$ depends only on the jump of $g\in H(F')$ across $F$, and $\Phi[g]$ can be interpreted in a manner more familiar to Numerical Analysts, as a \emph{hyperfunction} \cite{sato1959theory}. 

 In \cite{grothendieck1953certains}, Grothendieck generalised Silva's duality theorems to general complementary sets $F$ and $F'$, and to weakly holomorphic functions with values in a locally convex vector space $V$. Specifically, $H(F'; V) \cong H'(F;V')$, where $V'$ is the continuous dual space of $V$. We won't go into the details here about the strict definition of $H(F';V)$ --- we proved Theorem \ref{thm:main} without appealing to Grothendieck's theory --- because we merely need to know that the kernel $K$ in Theorem \ref{thm:main} lies in $ H(F'; \, L_\lambda^2(D))$ (after subtracting off $K(x,\infty)$ if $\infty \in F'$, which does not change $\mathcal{K}'$), and
\begin{equation}
    H(F'; \, L_\lambda^2(D)) \cong H'(F; \, L_\lambda^2(D)).
\end{equation}
From these two facts, we deduce the existence of a corresponding distributional kernel $K' \in H'(F; \, L_\lambda^2(D))$. Specifically, there exists a distribution on analytic functions $K'$ such that
\begin{equation}\label{eqn:KandK'}
    \frac{1}{2\pi \mathrm{i}}\int_\Gamma K(x,\xi) \, h(\xi) \, \mathrm{d}\xi = \int_F K'(x,z) \, h(z) \, \mathrm{d}\nu(z),
\end{equation}
for all $h \in H(F)$. There are two interpretations of the space $H'(F; \, L_\lambda^2(D))$ in this equation, one through Grothendieck duality and one through the Gelfand triple (or rigged Hilbert space),
\begin{equation}
    H(F; \, L_\lambda^2(D)) \subset L_\nu^2(F; \, L_\lambda^2(D)) \subset H'(F; \, L_\lambda^2(D)),
\end{equation}
each with a corresponding kernel. Setting $h(z) = (y-z)^{-1}$ for a fixed $y \in F'$ yields the key relationship,
\begin{equation}
    K(x,y) = \int_F \frac{K'(x,z)}{z-y} \, \mathrm{d} \nu(z).
\end{equation}
If the distribution $K'$ is actually a function, then we can estimate the norm of $\mathcal{K}'$ as follows.

\begin{lemma}\label{lem:K'}
If $K'\in L^{p}_\nu(F; \, L_\lambda^q(D))$, where $p,q \in [1,\infty]$, then
    \begin{equation*}
        \|\mathcal{K}'\|_{H_\nu^{p^*}(F) \to L_\lambda^q(D)} \leq \|K'\|_{L_\nu^p(F; \, L_\lambda^q(D))},
    \end{equation*}
    where $p^*$ is the H\"older dual exponent, satisfying $\tfrac{1}{p} + \tfrac{1}{p^*} = 1$.
    \begin{proof}
This follows from Equation \eqref{eqn:KandK'} and Lemma \ref{lem:boundZ}.
    \end{proof}
\end{lemma}

\section{Computation by rational interpolation}\label{sec:rational}
The low-rank approximant used in the proof of Theorem \ref{thm:main} is a Hilbert--Schmidt operator with kernel
\begin{equation}\label{eqn:Kn}
            K_n(x,y) = \frac{1}{2\pi \mathrm{i}} \int_\Gamma K(x,\xi) \, \left(1 - \frac{\phi(y)}{\phi(\xi)}\right) \frac{1}{y-\xi} \, \mathrm{d} \xi,
        \end{equation}
where $\phi(\xi) = \prod_{i=1}^\ell (\xi-q_i) / \prod_{k=1}^m (\xi - p_k)$ for some $\ell,m \leq n$. This formula bears a resemblance to the Hermite integral formula for rational interpolation with prescribed poles.
\begin{theorem}[Hermite integral formula {\cite[Thm 8.2]{walsh1935interpolation}}]\label{thm:HIF}
	Let $\Omega$ be a bounded open subset of $\mathbb{C}$ whose boundary $\Gamma$ is a finite sum of rectifiable Jordan curves. Let $g \in H(\Omega)$, $q_1,\ldots, q_\ell \in \Omega$ and $p_1,\ldots,p_m \in \mathbb{C}\setminus \{q_1,\ldots,q_\ell\}$. There exists a unique $r \in \mathcal{R}_{\ell-1,m}$ that interpolates $g$ at $\{q_1,\ldots,q_\ell\}$ (in the Hermite sense if points are repeated) and has poles at $\{p_1,\ldots,p_m\}$, which for $y \in \Omega$ satisfies
	\begin{equation*}
		r(y) = \frac{1}{2\pi \mathrm{i}} \int_\Gamma g(\xi) \,  \left(1 - \frac{\phi(y)}{\phi(\xi)}\right) \frac{1}{\xi-y} \, \mathrm{d} \xi.
	\end{equation*}
\end{theorem}
Theorem \ref{thm:HIF} indeed implies, after changing the orientation of the contour and reversing the sign of the $\xi - y$ term, that the low-rank kernel given in the proof of Theorem \ref{thm:main} is a rational interpolant of $K(x,\cdot)$, with interpolation nodes and poles that are independent of $x$. The roots and poles in the proof are chosen so that $\phi$ is optimal for $Z_n(L_\mu^2(E),L_\nu^p(F))$. If the roots of $\phi$ are distinct, then
    \begin{equation}\label{eqn:Knexpansion}
        K_n(x,y) = \sum_{j=1}^\ell K(x,q_j) \, \phi_j(y), \text{ where } \phi_j(y) =  \prod_{\substack{i=1 \\ i \neq j}}^\ell \frac{y-q_i}{q_j-q_i} \, \prod_{k=1}^m \frac{q_j-p_k}{y-p_k}.
    \end{equation}
    To evaluate $K_n$ in an efficient and numerically stable manner, use the barycentric form \cite{higham2004numerical, webb2012stability},
\begin{equation}\label{eqn:ratbary}
    K_n(x,y) = \phi(y)\,\sum_{j=1}^\ell \frac{K(x,q_j)\, w_j }{y-q_j}, \text{ where } w_j = \frac{\prod_{k=1}^m (q_j-p_k)}{\prod_{\substack{i=1 \\ i \neq j}}^\ell (q_j - q_i)}.
\end{equation}
The optimality of both the interpolation nodes and the poles is crucial, and is a key aspect distinguishing this approach from other rational approximations such as the proxy point method \cite[Eqn~(2.3)]{ye2020analytical}. Now the only thing standing between us and a practical algorithm is the choice of $\phi$. The case of two disjoint intervals on the extended real line has an explicit solution due to Zolotarev.
\begin{theorem}[Zolotarev \cite{zolotarev1877application}, {\cite[p.~20]{wilber2021computing}}]\label{thm:Zolotarevphi}
Let $\phi(z)$ be the Zolotarev rational function for
$Z_n([a,b],[c,d])$, where $-\infty \leq a < b < c < d \leq \infty$. Then the zeros $q_1,\ldots,q_n$ and poles $p_1,\ldots,p_n$ of $\phi$ are
\begin{equation}
q_j = T \left(-\mathrm{dn} \left[ \frac{2j - 1}{2n} \mathrm{K}(\Xi), \Xi \right] \right), \quad
p_j = T\left(\mathrm{dn} \left[ \frac{2j - 1}{2n} \mathrm{K}(\Xi), \Xi \right] \right),
\end{equation}
where $\Xi = \sqrt{1 - 1/\tau^2}$, $\tau = -1 + 2 \gamma + 2 \sqrt{\gamma^2 - \gamma}$, $\gamma$ is as in Theorem \ref{thm:BT}, $\mathrm{K}$ is the complete elliptic integral of the first kind \cite[19.2.8]{NIST:DLMF}, $\mathrm{dn}(z, \Xi)$ is the Jacobi elliptic function of the third kind \cite[22.2.6]{NIST:DLMF}, and
$T$ is the M\"obius transformation such that
$T(-1) = a, T(-1/\tau) = b, T(1/\tau) = c, T(1) = d$.
\end{theorem}
In this case, $q_1,\ldots,q_n$ are distinct, so formula \eqref{eqn:ratbary} can be used to evaluate the kernel $K_n$. However, for other sets $E$ and $F$, the roots of the optimal $\phi$ may not be distinct. In such cases, the rational interpolant is of Hermite type and requires partial derivatives of $K$ with respect to $y$. Indeed, for $R > 1$, the optimal Zolotarev rational function for $Z_n(\{|z| \leq 1 \}, \{|z| \geq R \})$ is $\phi(\xi) = \xi^n$. 

However, the rational function $\phi$ with $q_j =  \mathrm{e}^{2\pi\mathrm{i}j/n}$ and $p_j = R \mathrm{e}^{2\pi\mathrm{i}j/n}$ is asymptotically optimal, in the sense that $\lim_{n \to \infty} \sup_{y\in E, z \in F} |\phi(z)^{-1}\phi(y)|^{1/n} = \rho(E,F)^{-1}$ \cite{walsh1965hyperbolic, starke1992near}. Asymptotically optimal rational functions with distinct roots can be computed using conformal mapping to an annulus (if possible) \cite{walsh1965hyperbolic, starke1992near}, Bagby's iteration \cite{bagby1969interpolation}, or the recent algorithm of Trefethen and Wilber \cite{trefethen2025computation} based on the AAA-Lawson algorithm applied to the sign function on $E \cup F$.

Of course, suboptimal $\phi$ leads to suboptimal interpolants. The following suboptimal version of Theorem \ref{thm:main} places the bounds of Little--Reade type (based on polynomial approximation) in the context of the present work \cite{little1984eigenvalues} (see also \cite[Thm 7.1]{townsend2015continuous}).
\begin{proposition}\label{prop:LR}
Let $K \in C(D\times [-1,1])$ be as in Theorem \ref{thm:main}, such that $K(x,\cdot)$ is analytically continuable to the open interior of a Bernstein ellipse with parameter $R > 1$. Then for $n = 1,2,\ldots$, there exists a rank-$n$ kernel $K_n \in C(D\times E)$ such that, for any $p \in [1,\infty]$,
    \begin{equation*}
    \|\mathcal{K}-\mathcal{K}_n\|_{L_\mu^2(-1,1) \to L_\lambda^2(D)} \leq \frac{2R^{-n}}{1-R^{-2n}} \, \left\| \mathcal{C}\right\|_{L_\mu^2(-1,1) \to L_\nu^p(F)} \,\|\mathcal{K}'\|_{H_\nu^p(F) \to L_\lambda^2(D)},
\end{equation*}
where $\mathcal{C}$ is the Cauchy transform in Equation \eqref{eqn:Cauchytransform}.
\begin{proof}
     Apply Theorem \ref{thm:main} and take $\phi(\xi) = T_n(\xi) \in \mathcal{R}_{n,0} \subset \mathcal{R}_{n}$, which corresponds to polynomial interpolation at the roots of $T_n$. To bound the resulting suboptimal Zolotarev term, start with the formula $\phi(z) = \tfrac12(\zeta^{n}+\zeta^{-n})$ where $z = \tfrac12(\zeta + \zeta^{-1})$ with $|\zeta| \geq 1$ \cite[18.5.1]{NIST:DLMF}. Then for $y \in [-1,1]$, $z \in F$,
     \begin{equation}
         \frac{\phi(z)^{-1} \phi(y)}{y-z} = \frac{2\zeta^{-n}}{1+\zeta^{-2n}} \cdot C(z,y) \cdot T_n(y), \quad \text{(where } C \text{ is the kernel of } \mathcal{C}).
     \end{equation}
     To conclude, note that $|\zeta| \geq R$ for all $z \in F$ since $F$ lies outside the ellipse.
\end{proof}
\end{proposition}
\begin{remark}\label{rem:Cheb}
If $F \subset \mathbb{R} \cup \{\infty\}$ then the denominator can be changed to $1+R^{-2n}$ because $|T_n(z)| \geq \tfrac12(R^n + R^{-n})$ for $z \in F$. If $\mathcal{C}$ is unbounded then $\phi(\xi) = (\xi-t)T_{n-1}(\xi)$ can be used as in Lemma \ref{lem:Z1}.
\end{remark}
\section{Examples}\label{sec:examples}

MATLAB scripts to generate all of the figures in this paper are available at the GitHub repository, \href{https://github.com/marcusdavidwebb/LRAAK}{https://github.com/marcusdavidwebb/LRAAK}.

\subsection{Cauchy matrices and tensors}

Consider a Cauchy matrix $A \in \mathbb{R}^{N\times N}$
\begin{equation}\label{eqn:Cauchyexample1}
    A_{i,j} = \frac{1}{x_i + y_j}, \qquad x_i \in D, \quad y_j \in E,
\end{equation}
where $D$ and $E$ are compact subsets of the complex plane such that $(-D) \cap E = \emptyset$. An analysis of low-rank approximations to this kernel (including lower bounds) appears in \cite{oseledets2007lower}. To apply Theorem \ref{thm:main}, we take the kernel $K(x,y) = \frac{1}{x+y}$, with $\lambda$ and $\mu$ being discrete counting measures on $\{x_i\}$ and $\{y_j\}$. For each $x \in D$, $K(x,\cdot)$ can be analytically continued to $\Chat \setminus \{-x\}$, so we can take $F = -D$. The Grothendieck dual of $\mathcal{K}$ is $\mathcal{K}':H(-D) \to C(D)$, given by
\begin{equation}
    \mathcal{K}'[h](x_i) = \frac{1}{2\pi\mathrm{i}}\int_\Gamma \frac{h(\xi)}{x_i + \xi} \, \mathrm{d}\xi = h(-x_i),
\end{equation}
where $\Gamma$ is any suitable contour within the region of analyticity of $h$, winding once around $-D$ (and zero times around $E)$. It is clear that $\mathcal{K}'$ merely reflects $h$ across the origin, so if we take $\mathrm{d}\nu(z) = \mathrm{d}\lambda(-z)$, then
$\|\mathcal{K}'\|_{H_\nu^2(-D) \to L^2_\lambda(D)} = 1$. Therefore, by Theorem \ref{thm:main},
\begin{eqnarray}
    \sigma_{n+1}(A) \leq \|\mathcal{K}-\mathcal{K}_n\|_{L^2_\mu(E)\to L^2_\lambda(D)} &\leq& Z_n(L_\mu^2(E),L_\nu^2(-D)) \\
    &\leq& Z_n(E,-D) \, \|(y-z)^{-1}\|_{L_\mu^2(E)\to L_\nu^2(-D)},\label{eqn:Cauchy1}
\end{eqnarray}
by Lemma \ref{lem:boundZ}. Since $\|(y-z)^{-1}\|_{L_\mu^2(E)\to L_\nu^2(-D)} = \|(y+z)^{-1}\|_{L_\mu^2(E)\to L_\lambda^2(D)} =  \sigma_1(A)$, we recover the Beckermann--Townsend results exactly for Cauchy matrices \cite{beckermann2019bounds}.

Now consider a Cauchy tensor $A \in \mathbb{R}^{N\times N \times N}$,
\begin{equation}\label{eqn:Cauchexample2}
    A_{i,j,k} = \frac{1}{w_i + x_j + y_k}, \qquad (w_i,x_j) \in D, \quad y_k \in E,
\end{equation}
where $E$ is a compact subset of $\mathbb{C}$ and $D$ is a compact subset of $\mathbb{C}^2$ such that $(-\Sigma D) \cap E = \emptyset$, where $\Sigma D =\{ w+x : (w,x) \in D\}$. Then we take the kernel $K((w,x),y) = \frac{1}{w + x + y}$, with $\lambda$ and $\mu$ counting measures as above. For each $(w,x) \in D$, $K((w,x),\cdot)$ can be analytically continued to $\Chat\setminus \{-(w+x)\}$. Therefore, we take $F = -\Sigma D$. The Grothendieck dual of $\mathcal{K}$ is $\mathcal{K}' : H(-\Sigma D) \to C(D)$ and we have
\begin{equation}
    \mathcal{K}'[h]((w,x)) = \frac{1}{2\pi \mathrm{i}} \int_\Gamma \frac{h(\xi)}{w + x + \xi} \, \mathrm{d} \xi = h(-w - x)
\end{equation}
The norm is slightly more challenging to bound. If we define a Radon measure $\nu$ on $-\Sigma D$ by
\begin{equation}
    \nu(U) = \mu\left( \{ (w,x) \in D : -w-x \in U \} \right),
\end{equation}
for all open sets $U \subset -\Sigma D$, then $\mathcal{K}' : H_\nu^2(-\Sigma D) \to L_\lambda^2(D)$ is an isometry. Therefore, by Theorem \ref{thm:main},
\begin{equation}
    \|\mathcal{K}-\mathcal{K}_n\|_{L^2_\mu(E)\to L^2_\lambda(D)} \leq Z_n(L_\mu^2(E),L_\nu^2(-\Sigma D)) \leq Z_n(E, -\Sigma D)\, \|(y-z)^{-1}\|_{L_\mu^2(E)\to L_\nu^2(-\Sigma D)}.\label{eqn:Cauchy2}
\end{equation}
Just as with the Cauchy matrices, $\|(y-z)^{-1}\|_{L_\mu^2(E)\to L_\nu^2(-\Sigma D)} = \|\mathcal{K}\|_{L_\mu^2(E) \to L_\lambda^2(D)}$, so we recover the results of Shi and Townsend \cite{shi2021compressibility} based on displacement structure of the Cauchy tensor.

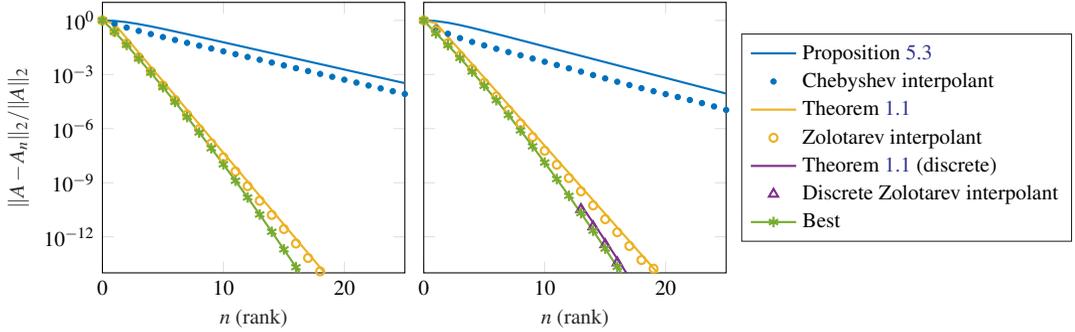
\begin{figure}[!t]
\centering
\footnotesize
\input{figs/figure3.tex}
\caption{Cauchy matrices and tensors. On the left I used $A\in\mathbb{R}^{N\times N}$ in Equation \eqref{eqn:Cauchyexample1} with $N = 100$, $D = [1, 70]$, $E = [2,100]$. On the right I used $A \in \mathbb{R}^{N\times N \times N}$ in Equation \eqref{eqn:Cauchexample2} with $N=50$, $D = [1,70]\times[1,199] \subset \mathbb{C}^2$, $E = [2,100]$. In both plots the sample points are equispaced over their respective intervals.}\label{fig:cauchyfigs}
\end{figure}

Figure \ref{fig:cauchyfigs} shows the results of numerical experiments for a Cauchy matrix and a Cauchy tensor. The blue dots are the errors associated with the polynomial interpolants to $K(x,\cdot)$ at Chebyshev nodes on $[2,100]$ and the blue lines are application of Proposition \ref{prop:LR} mapped to $[2,100]$. The yellow circles are the errors associated with rational interpolants with nodes and poles associated with $Z_n([2,100],[-100,-2])$ as in Theorem \ref{thm:Zolotarevphi} and the yellow lines correspond to Equation \eqref{eqn:Cauchy1} and \eqref{eqn:Cauchy2} with $Z_n(E,F)$ bounded by Theorem \ref{thm:BT}. The purple triangles are the errors associated with the Zolotarev rational functions for the discrete sets $E = \{y_j\}$ and $F = \{-{w_k - x_i}\}$, rather than on the continuous intervals, and yielded an improved approximation, almost indistinguishable from the best. These discrete Zolorev rational functions were computed using \cite{trefethen2025computation}, which also returns the associated Zolotarev numbers $Z_n(E,F)$ yielding the purple lines when substituted into \eqref{eqn:Cauchy1} and \eqref{eqn:Cauchy2}.

\subsection{Log-Cauchy matrices}\label{subsec:logCauchy}
Consider a matrix $A \in \mathbb{R}^{N\times N}$ given by
\begin{equation}\label{eqn:logCauchy}
A_{i,j} = \log(x_i + y_j),
\end{equation}
where $x_1,\ldots, x_N, y_1,\ldots,y_N \in [c,d]$, where $d > c > 0$. This matrix samples the kernel
\begin{equation}
    K(x,y) = \log(x+y) \qquad x,y\in [c,d],
\end{equation}
where $0 < c < d < \infty$. For all $x \in [c,d]$, $K(x,\cdot)$ is analytic in $\Chat \setminus [-\infty, -x]$, so we can take $D = E = [c,d]$, $F = [-\infty, -c]$. Let us examine the Grothendieck dual operator. If $h \in H^1([-\infty,-c])$, then following the same proof as in Lemma \ref{lem:logkernel} (except that $a_1 = 0$ because $h \in L^1(-\infty,-c)$), we have
\begin{equation}\label{eqn:logK'}
 \mathcal{K}'[h](x) = - \int_{-\infty}^{-x} h(z) \, \mathrm{d} z.
 \end{equation}
The kernel of the Grothendieck dual is therefore $K' \in L^1([-\infty,-c];\, L_\lambda^2([c,d]))$, with
    \begin{equation}
        K'(x,z) = \begin{cases}
        -1 & \text{ if } x+z \leq 0 \\
        0 & \text{ otherwise,}
    \end{cases}
    \end{equation}
where $x \in [c,d]$, $z \in [-\infty,-c]$. By Lemma \ref{lem:K'} with $p = \infty$, $q = 2$, 
\begin{equation}
    \|\mathcal{K}'\|_{H^1([-\infty,-c]) \to L_\mu^2(c,d)} = \sup_{z < -c} \left(\int_c^d |K'(x,z)|^2 \,\mathrm{d}\lambda(x)\right)^{1/2} = \lambda([c,d])^{1/2}.
\end{equation}
By Lemma \ref{lem:Z1}, we have
\begin{equation}
    Z_n(L_\mu^2(c,d), L^1(-\infty,-c)) \leq \frac12\log\left(\frac{d+c}{2c} \right) \, \mu([c,d])^{1/2} \,  Z_{n-1}([c,d], [-\infty, -c]).
\end{equation}
Therefore, by Theorem \ref{thm:main}, in the case that $\lambda$ and $\mu$ are counting measures,
\begin{equation}\label{eqn:logCauchybound}
    \|\mathcal{K} - \mathcal{K}_n \|_{L_\mu^2([c,d]) \to L_\mu^2([c,d])} \leq \frac{N}{2}\,\log\left(\frac{d+c}{2c} \right) \,  Z_{n-1}([c,d], [-\infty, -c]).
\end{equation}

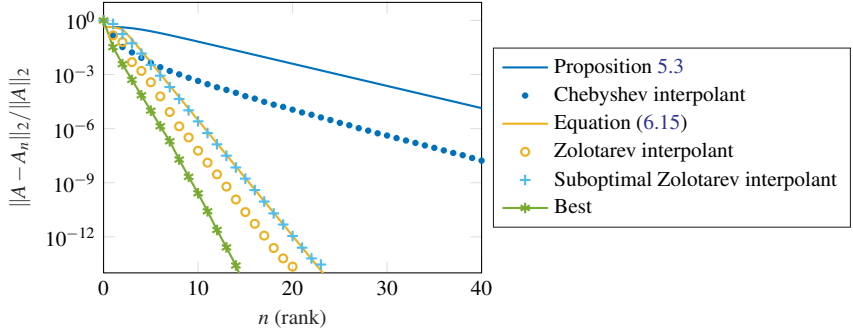
\begin{figure}[!t]
\centering
\footnotesize
\input{figs/figure4.tex}
\caption{Log-Cauchy matrices. I used $A \in \mathbb{R}^{N\times N}$ in Equation \eqref{eqn:logCauchy} with $N = 100$, $D = E = [1,N]$ and $\{x_i\}$, $\{y_j\}$ randomly uniformly distributed in $[1,N]$ (and $x_1 = y_1 = 1$, $x_N = y_N = N$). The Chebyshev interpolant (blue dots) uses the modified interpolation points mentioned in Remark \ref{rem:Cheb}. The Zolotarev interpolant (yellow circles) uses the $\phi$ in the proof of Lemma \ref{lem:Z1} for the nodes and poles, whereas the suboptimal Zolotarev interpolants (cyan plus signs) uses the $\phi$ that is optimal for $Z_n(E,F)$.}\label{fig:logkernel}
\end{figure}

Figure \ref{fig:logkernel} shows the results of numerical experiments for a log-Cauchy matrix with randomly sampled $\{x_i\}$ and $\{y_j\}$ in $[c,d] = [1,N]$. The blue dots and lines are as in Figure \ref{fig:cauchyfigs}. The yellow circles are the approximation error associated with rational interpolants with nodes and poles discussed in the proof of Lemma \ref{lem:Z1}, that is, instead of the optimal nodes and poles for $Z_n(L_\mu^2(c,d), L^1(-\infty,-c))$, we take the optimal nodes and poles associated with $Z_{n-1}([c,d],[-\infty,-c])$ with an extra node at $-c + \sqrt{(d+c)\cdot 2c}$. The cyan plus signs are the errors associated with the rational interpolant at the nodes and poles associated with $Z_n([c,d],[-\infty,-c])$, which are about 50 times worse than the yellow circles.

\subsection{Twisted Hankel transform matrix}

The discrete Hankel transform is defined by the matrix \cite{johnson1987improved},
\begin{equation}
    J_{i,j} = \mathrm{J}_\nu\left( \frac{\omega_i\,\omega_j}{\omega_{N+1}} \right), \qquad i,j = 1,2,\ldots,N,
    \end{equation}
where $\omega_1,\ldots,\omega_{N+1}$ are the first $N+1$ positive roots of $\mathrm{J}_\nu$, the Bessel function of order $\nu \geq 0$ \cite[§10]{NIST:DLMF}. For simplicity, we take $\nu = 0$. A new class of fast algorithms for the discrete Hankel transform (to appear in a forthcoming publication) rely on the low-rank approximation of the kernel
\begin{equation}\label{eqn:HankelKernel}
    K(x,y) = \mathrm{H}_0^{(1)}(xy)\mathrm{e}^{-\mathrm{i}xy}, \qquad x, y \in \left[\widehat{\omega}_1, \widehat{\omega}_N\right], \qquad \widehat{\omega}_j = \omega_j/\sqrt{\omega_{N+1}},
\end{equation}
where $\mathrm{H}_0^{(1)}$ is the Hankel function of the first kind. For each $x \in [\widehat{\omega}_1,\widehat{\omega}_N]$, $K(x,\cdot)$ is analytic in the set $\mathbb{C} \setminus (-\infty,0]$, where $\mathrm{H}_0^{(1)}$ has its principle branch cut. Therefore, we can take $D = E = [\widehat{\omega}_1,\widehat{\omega}_N]$ and $F = [-\infty, 0]$. Its Grothendieck dual $\mathcal{K}':H([-\infty,0]) \to C([\widehat{\omega}_1,\widehat{\omega}_N])$ is
    \begin{equation}
    \mathcal{K}'[h](x) = \frac{2}{\pi \mathrm{i}} \int_{-\infty}^0 \mathrm{J}_0(xz)\, \mathrm{e}^{-\mathrm{i}xz} \,h(z) \, \mathrm{d} z.
    \end{equation}
The proof of this is similar to that of Equation \eqref{eqn:logK'}. The two consequential things to note are, that $\mathrm{H}_0^{(1)}(\xi)$ has a jump of $-4\mathrm{J}_0(\xi)$ over the negative real axis \cite[10.11.5]{NIST:DLMF}, and, crucially,
        \begin{equation}
            \mathrm{H}_0^{(1)}(\xi)\mathrm{e}^{-\mathrm{i}\xi} = \mathcal{O}(\xi^{-1/2}) \quad \text{ as } \xi \to \infty \text{ with } |\mathrm{arg}(\xi)| < \pi-\delta,
        \end{equation}
        by \cite[10.7.8]{NIST:DLMF}. This second fact ensures that the contour integrals depending on $R$ converge to zero as $R \to \infty$. Therefore, as a member of $H'([-\infty,0]; L_\lambda^2([\widehat{\omega}_1,\widehat{\omega}_N]))$, the Grothendieck dual has kernel $K'(x,z) = \frac{2}{\pi\mathrm{i}} \mathrm{J}_0(xz) \mathrm{e}^{-\mathrm{i}xz}$, so we can use Lemma \ref{lem:K'} to bound its norm:
\begin{equation}
    \|\mathcal{K}'\|_{H^1([-\infty,0]) \to L_\lambda^2([\widehat{\omega}_1,\widehat{\omega}_N])} = \sup_{z \leq 0} \frac{2}{\pi}\left( \int_{\widehat{\omega}_1}^{\widehat{\omega}_N} |\mathrm{J}_0(xz)|^2 \, \mathrm{d}\lambda(x)\right)^{1/2} = \frac{2}{\pi} \lambda([\widehat{\omega}_1,\widehat{\omega}_N])^{1/2},
\end{equation}
because $\mathrm{J}_0(0) = 1$. By Lemma \ref{lem:Z1}, we have
\begin{equation}
    Z_n(L_\mu^2([\widehat{\omega}_1,\widehat{\omega}_N]), L^1(-\infty,0)) \leq \frac12\,\log\left(\frac{\omega_N}{\omega_1} \right) \, \mu([\widehat{\omega}_1,\widehat{\omega}_N])^{1/2} \,  Z_{n-1}\left([-\infty, 0], \left[1,\frac{\omega_N}{\omega_1}\right] \right).
\end{equation}
Therefore, by Theorem \ref{thm:main}, in the case that $\lambda$ and $\mu$ are counting measures,
\begin{equation}\label{eqn:Hankelbound}
     \|\mathcal{K} - \mathcal{K}_n \|_{L_\mu^2([\widehat{\omega}_1,\widehat{\omega}_N]) \to L_\lambda^2([\widehat{\omega}_1,\widehat{\omega}_N])} \leq \frac{N}{\pi}\,\log\left(\frac{\omega_N}{\omega_1} \right) \,  Z_{n-1}\left([-\infty, 0], \left[1,\frac{\omega_N}{\omega_1}\right] \right).
\end{equation}

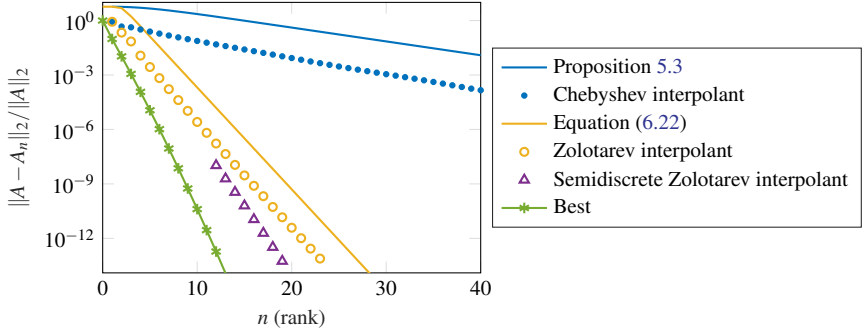
\begin{figure}[!t]
\centering
\footnotesize
\input{figs/figure5.tex}
\caption{Twisted Hankel transform matrix. I used $A \in \mathbb{R}^{N\times N}$ in Equation \eqref{eqn:twistedHankelmatrix} with $N = 100$. The Chebyshev interpolant (blue dots) uses the modified interpolation points mentioned in Remark \ref{rem:Cheb}. The Zolotarev interpolant (yellow circles) uses the $\phi$ in the proof of Lemma \ref{lem:Z1} for the nodes and poles, whereas the semidiscrete Zolotarev interpolants (purple triangles) use the $\phi$ that is optimal for $Z_n(\{\widehat{\omega}_1,\widehat{\omega}_2,\ldots,\widehat{\omega}_N\},[-\infty,0])$ as described in the main text.}\label{fig:hankelkernel}
\end{figure}

Figure \ref{fig:hankelkernel} shows the results of numerical experiments for the matrix, \begin{equation}\label{eqn:twistedHankelmatrix}
    A_{i,j} = K(\widehat{\omega}_i, \, \widehat{\omega}_j) \qquad i,j = 1,2,\ldots, N.
\end{equation}
The blue dots, blue line, yellow circles and yellow line are as in Figure \ref{fig:logkernel}. The purple triangles correspond to the interpolation points and poles taken to be the roots and poles of the Zolotarev rational function for $Z_{n}(\{\widehat{\omega}_1,\widehat{\omega}_2,\ldots,\widehat{\omega}_N\},[-\infty,0])$. To compute this $\phi$, I first mapped these sets to bounded sets with the M\"obius transformation $T(z) = \frac{1}{1 + \widehat{\omega}_N-z}$. I used Trefethen and Wilber's algorithm \cite{trefethen2025computation}, which requires discrete sets as input, so I approximated $T([-\infty, 0])$ by 2000 equally spaced points between $0$ and $(1+ \widehat{\omega}_N)^{-1}$. Then I mapped the resulting zeros and poles (which are covariant under M\"obius transformations) back to the original sets. The resulting interpolants (purple triangles) achieved a smaller error than those implied by Equation \eqref{eqn:Hankelbound} (yellow circles).

For the semidiscrete interpolant, I did not use the extra interpolation point $t$ from the proof of Lemma \ref{lem:Z1}, because the results were, for an unknown reason, slightly worse. The results of this paper belie some difficulties in obtaining consistent results from Trefethen and Wilber's algorithm (which the authors describe as ``90\% reliable'' \cite{trefethen2025computation}).

\subsection{Beta-Cauchy matrices}\label{sec:BetaMatrices}
Positive-definite Hankel matrices of the following form appear in the change-of-basis matrices between classical orthogonal polynomials \cite{townsend2018fast}. A family of fast algorithms for performing the matrix-vector product depends on low-rank approximation of the matrix,
\begin{equation}
    A_{ij} = \mathrm{B}(i+j+\alpha,\beta), \qquad  i,j=0,\ldots,N,
\end{equation}
where $\alpha,\beta > 0$, $\beta \notin\mathbb{N}$, and $\mathrm{B}$ is Euler's beta function, $\mathrm{B}(\alpha,\beta) = \Gamma(\alpha)\Gamma(\beta)/\Gamma(\alpha+\beta)$ \cite[5.12]{NIST:DLMF}. To apply Theorem \ref{thm:main}, we take
\begin{equation}
    K(x,y) = \mathrm{B}(x+y+\alpha,\beta), \qquad x,y \in [0,N],
\end{equation}
with $\lambda$ and $\mu$ the counting measures on $\{0,1,\ldots, N\}$. For each $x \in \{0,\ldots,N\}$, $\xi \mapsto K(x,\xi)$ is analytic for $\xi \in \mathbb{C} \setminus \{  - (\alpha + x + k): k \in \mathbb{Z}_{\leq 0} \}$, so we can take
$D = [0,N]$, $E = [0,N]$, $F = [-\infty, -\alpha]$. We use the partial fraction expansion \cite[p.~8]{bateman1953higher},
\begin{equation}
    \mathrm{B}(\alpha,\beta) = \frac{1}{\Gamma(1-\beta)}\,\sum_{k=0}^\infty \frac{\Gamma(k+1-\beta)}{\Gamma(k+1)} \,\frac{1}{\alpha+k},
\end{equation}
to obtain
\begin{equation}
    K(x,y) = \sum_{k=0}^\infty  \frac{w(k,\beta)}{x+y+\alpha+k}, \qquad w(k,\beta) =\frac{\Gamma(k+1-\beta)}{\Gamma(k+1)\,\Gamma(1-\beta)}.
\end{equation}
For intuition in what follows, it is useful to note by Stirling's formula that $w(k,\beta) \sim k^{-\beta}/\Gamma(1-\beta)$ for large $k$. The Grothendieck dual operator is, for $h \in H([-\infty,-\alpha])$,
\begin{equation}
    \mathcal{K}'[h](x) = \sum_{k=0}^\infty w(k,\beta) \, h(-\alpha - k - x) = \sum_{k=x}^\infty w(k-x,\beta) \, h(-\alpha - k).
\end{equation}
An argument similar to that in Section \ref{sec:proofs} is technically required to preclude hiccoughs at infinity. If we choose the space $H_\nu^1([-\infty,-\alpha])$, where $\nu$ is the discrete weighted measure, $\nu(z) = \sum_{k=0}^\infty |w(k,\beta)| \, \delta (z + k + \alpha)$, then
\begin{equation}
    \mathcal{K}'[h](x) = \int_{-\infty}^{-\alpha} K'(x,z)\,h(z) \, \mathrm{d} \nu(z),
\end{equation}
where, for $x \in \{0,1,\ldots,N\}$ and $z \in \{-\alpha, -1-\alpha,-2-\alpha,\ldots\}$,
\begin{equation}
    K'(x,z) = \begin{cases}
         \frac{w(-z-\alpha-x,\beta)}{|w(-z-\alpha,\beta)|} & \text{ if }  x+z + \alpha \leq 0\\ 
         0 & \text{ otherwise.}
    \end{cases}
\end{equation}
The weighting of the discrete measure $\nu$ has provided a uniformly bounded $K'(x,z)$, since for $k \geq x$,
\begin{equation}
    |K'(x,-k-\alpha)| = \left|\frac{\Gamma(k+1-x-\beta)}{\Gamma(k+1-\beta)} \frac{\Gamma(k+1)}{\Gamma(k+1-x)}\right| = \prod_{\ell=0}^{x-1} \frac{k-\ell}{k-\ell - \beta} = \prod_{\ell=0}^{x-1}\frac{1}{1-\frac{\beta}{k-\ell}}.
\end{equation}
This expression is strictly decreasing with $k$, so using Lemma \ref{lem:K'} we can bound the norm of $\mathcal{K}'$ by
\begin{eqnarray}
    \|\mathcal{K}'\|_{H^1_\nu([-\infty,-\alpha]) \to L_\lambda^2([0,N])} &\leq& \sup_{z \in\{-\alpha,-1-\alpha,\ldots\}} \left(\sum_{x=0}^N |K'(x,z)|^2 \right)^{1/2} \\
    &=& \left(\sum_{x=0}^{N} |K'(x,-N-\alpha)|^2\right)^{1/2} \\
    &=& \left(\sum_{x=0}^{N} \left(\frac{w(x,\beta)}{w(N,\beta)}\right)^2\right)^{1/2}.
\end{eqnarray}
The Cauchy transform $\mathcal{C} : L_\mu^2([1,N]) \to L_\nu^1([-\infty,-\alpha])$ can also be bounded. By Lemma \ref{lem:boundZ},
\begin{eqnarray}
    \|\mathcal{C}\| &\leq& \left( \int_0^N \left(\int_{-\infty}^{-\alpha} |y-z|^{-1} \, \mathrm{d}\nu(z)\right)^2 \, \mathrm{d}\mu(y) \right)^{1/2} = \left( \sum_{y=0}^N \left( \sum_{k=0}^\infty \frac{|w(k,\beta)|}{y+k+\alpha} \right)^2 \right)^{1/2}.
\end{eqnarray}
Theorem \ref{thm:main} and Lemma \ref{lem:ZntoZn-1} imply
\begin{equation}\label{eqn:Betabound}
    \sigma_{n+1}(A) \leq Z_n([-\infty,-\alpha],[0,N]) \underbrace{\left( \sum_{y=0}^N \left( \sum_{k=0}^\infty \frac{|w(k,\beta)|}{y+k+\alpha} \right)^2 \right)^{1/2}\left(\sum_{x=0}^{N} \left(\frac{w(x,\beta)}{w(N,\beta)}\right)^2\right)^{1/2}}_{\mathcal{O}(N)}.
\end{equation}
The underbraced term is $\mathcal{O}(N)$ since $w(x,\beta) \sim x^{-\beta}/\Gamma(1-\beta)$ for large $x$, by Stirling's formula. In the Introduction, Figure \ref{fig:Hankelintro} shows the results of numerical experiments in the case $N = 100$, $\alpha = 1/2$, $\beta = 1/2$. The yellow line is Equation \eqref{eqn:Betabound}. The purple line is Equation \eqref{eqn:Betabound} with $Z_n([-\infty,-\alpha],[0,N])$ replaced by a numerical approximation of $Z_n(\{-\tfrac12, -\tfrac32, -\tfrac52, \ldots, -\infty\}, \{0,1,\ldots, N\})$. The interpolation nodes and poles for the interpolants corresponding to the purple triangles were computed in a similar manner to those in Figure \ref{fig:hankelkernel}, using a M\"obius transformation and Trefethen and Wilber's algorithm \cite{trefethen2025computation}.

\section{Conclusion}\label{sec:conclusion}

This paper provides a new theorem to bound the best low-rank approximation error of a matrix or operator whose kernel is continuous and analytically continuable in one of its variables to a larger region of the complex plane. Key concepts introduced are the Cauchy--Zolotarev numbers and the Grothendieck dual operator. We described practical techniques to bound the error, showed that rational interpolation computes low-rank approximations that achieve the bounds, and confirmed the theory with several numerical examples.

The core of the proof is the operator decomposition $\mathcal{K} = \mathcal{K}'\,\mathcal{C}$, where $\mathcal{C}:L_\mu^2(E)\to H(F)$ is a Cauchy transform operator and $\mathcal{K}' : H(F) \to C(D)$ is the Grothendieck dual operator of $\mathcal{K}$. This reveals a latent displacement structure in $\mathcal{K}$, since $\mathcal{C}$ satisfies the displacement equation,
\begin{equation}
    \mathcal{M}_F \,\mathcal{C} - \mathcal{C}\,\mathcal{M}_E = \mathcal{L},
\end{equation}
where $\mathcal{M}_E[g](y) = T(y)\,g(y)$, $\mathcal{M}_F[h](z) = T(z)\,h(z)$, and $\mathcal{L}[g](z) = -\sqrt{T'(z)}\, \int_E \sqrt{T'(y)} g(y) \, \mathrm{d}\mu(y)$ (a rank-1 operator), for any M\"obius transformation $T$ such that $T(E \cup F)$ is bounded, because
\begin{equation}
    \frac{T(z)-T(y)}{z-y} = \sqrt{T'(z)}\,\sqrt{T'(y)}.
\end{equation}
Rank-$n$ approximations $\mathcal{C}_n$ to $\mathcal{C}$, with kernels as in Lemma \ref{lemma:C_n}, as described in \cite{beckermann2019bounds, clouatre2025lifting}, yield rank-$n$ approximations $\mathcal{K}_n = \mathcal{K}'\mathcal{C}_n$ to $\mathcal{K}$. The kernel of $\mathcal{K}_n$ has explicit formulas given in Equation \eqref{eqn:Kn}, and \eqref{eqn:ratbary} in the case of distinct nodes $q_1,\ldots,q_\ell$.

The proof can easily be generalised to bound $\|\mathcal{K}-\mathcal{K}_n\|_{L_\mu^q(E) \to L_\lambda^r(D)}$ for any $q,r \in [1,\infty]$. The case $q = 1$, $r = \infty$ corresponds to the `max norm', a particularly tractable case. I have carefully avoided discussing higher-order singularities such as $K(x,y) = (y-x)^{-2}$. In this case $\mathcal{K}'[h](x) = h'(x)$, the norm of which cannot easily be bounded using the techniques of this paper (but Theorem \ref{thm:main2} applies nonetheless). A Sobolev space, $W_\nu^{p,s}(F)$, could be introduced, upon which $\mathcal{K}'$ is more obviously bounded, but then one is left to bound what could be called Cauchy--Sobolev--Zolotarev numbers, $Z_n(L_\mu^2(E),W_\nu^{p,s}(F))$.

The duality of $H(F)$ and $H(F')$ discussed in Section \ref{sec:Grothendieck} directly impacts Potential Theory and approximation of analytic functions. The standard potential-theoretic methods of bounding the integral in the Hermite integral formula for a function $f$ with singularities within the set $F$ lead to an error bound for each contour $\Gamma$, which for fixed $n$ often tends to infinity as $\Gamma$ approaches $\partial F$ \cite{levin2006potential}. Grothendieck duality takes $\Gamma$ out of the equation (after all, the equation is invariant under the choice of $\Gamma$), resulting in a single error bound involving a norm of $f$, interpreted as an element of $H'(F)$, multiplied by a Cauchy--Zolotarev number (or its polynomial analogue). This is the topic of a future paper.

The choice of interpolation points and poles is ripe for investigation. More research on computation of Zolotarev rational functions associated with continuous and, more importantly, discrete subsets of $\Chat$ is needed, since Trefethen and Wilber's algorithm is not 100\% reliable \cite{trefethen2025computation}. Figure \ref{fig:logkernel} highlights the potential of direct computation of the optimal rational functions for Cauchy--Zolotarev numbers. When $F$ is unknown, the best approach to finding good points and poles may involve the AAA algorithm \cite{nakatsukasa2018aaa} applied to $K$. Developing a practical interpolation scheme in this setting could have significant consequences, impacting any application employing low-rank approximation of analytic kernels.

\section*{Acknowledgements}
The author thanks Andrew Horning, Yanghong Huang, Daan Huybrechs, Jonas Latz, Yuji Nakatuskasa, Will Parnell, Endre S\"uli, Fran\c{c}oise Tisseur, Alex Townsend and Heather Wilber for helpful comments and discussion. Probing questions from Georg Maierhofer and Nick Trefethen led to improvements in the mathematics and the presentation. Many thanks also to the anonymous reviewers, who made excellent suggestions for improving the paper.

\bibliographystyle{plain}
\bibliography{refs}

\newpage

\appendix

\section{Proof of Lemma \ref{lem:logkernel}}\label{sec:proofs}
    \begin{proof}
    Fix $x \in [-N,-1]$. Since $h$ is analytic at and tends to zero at infinity (by the definition of $H([-\infty,-1])$), there exists a Laurent expansion $h(z) = \sum_{n=1}^\infty \frac{a_n(x)}{(z-x-1)^{n}}$,
    convergent for all $|z-x-1| > R$ for some $R > 0$. Set $0 < \varepsilon < x+1$ (sufficiently small) and consider the standard keyhole contour, $\Gamma = \Gamma_+ \cup \Gamma_- \cup \Gamma_R  \cup \Gamma_\varepsilon$, where $\Gamma_{\pm} = [-\sqrt{R^2-\varepsilon^2}+x,x] \pm \varepsilon \mathrm{i}$, $\Gamma_\varepsilon = \{x + \varepsilon \mathrm{e}^{\mathrm{i} \theta} : \theta \in [-\pi/2,\pi/2] \}$, $\Gamma_R = \{x + R \mathrm{e}^{\mathrm{i} \theta} : \theta \in [-\pi+\delta,\pi-\delta]\}$,
    where $\delta$ is set so that $\Gamma_R$ touches $\Gamma_\pm$ at its endpoints. 
        %Since $h \in H([-\infty, -1])$, we can ensure that $\epsilon$ is small enough so that $h$ is analytic in the interior  of $\Gamma$ (interior on the Riemann sphere). 
        
        By residue calculus, $\mathcal{K}'[(z-x-1)^{-1}](x) = - \log(1+x-x) = 0$, so we can assume that $a_1 = 0$. We calculate,
        \begin{eqnarray*}
            & & \frac{1}{2\pi\mathrm{i}}\int_{\Gamma_- \cup \Gamma_+} \log(\xi-x) h(\xi) \, \mathrm{d} \xi \\
            &=& \frac{1}{2\pi\mathrm{i}}\int_{-\sqrt{R^2-\varepsilon^2}+x}^{x} \log(z - \mathrm{i}\varepsilon - x)h(z-\mathrm{i}\varepsilon) - \log(z + \mathrm{i}\varepsilon-x)h(z+\mathrm{i}\varepsilon) \, \mathrm{d} z \\
            &\to& -\int_{-\infty}^{x} h(z) \, \mathrm{d} z \qquad \text{ as } \varepsilon \to 0 \text{ and then } R \to \infty,
        \end{eqnarray*}
        because the jump across the negative axis for the logarithm is $2\pi \mathrm{i}$. We further calculate,
        \begin{eqnarray*}
            \left|\frac{1}{2\pi\mathrm{i}}\int_{\Gamma_\varepsilon} \log( \xi-x) h(\xi) \, \mathrm{d} \xi\right| &\leq& \frac{\mathrm{len}(\Gamma_\varepsilon)}{2\pi} \sup_{\theta \in [-\pi/2,\pi/2]} |\log(\varepsilon \mathrm{e}^{\mathrm{i}\theta}) h(x + \varepsilon \mathrm{e}^{\mathrm{i}\theta})| \\
            &\leq& \frac{\varepsilon}{2} \sqrt{\log^2(\varepsilon) + \frac{\pi^2}{4} } \sup_{\theta \in [-\pi/2,\pi/2]} |h(x + \varepsilon \mathrm{e}^{\mathrm{i}\theta})|,
        \end{eqnarray*}
        which tends to zero as $\varepsilon \to 0$. Since $a_1 = 0$, we have that $h(\xi) = O(\xi^{-2})$ as  $|\xi| \to \infty$. Therefore,
        \begin{eqnarray*}
            \frac{1}{2\pi\mathrm{i}}\int_{\Gamma_R} \log(\xi-x) h(\xi) \, \mathrm{d} \xi &=& \frac{1}{2\pi}\int_{-\pi+\delta}^{\pi - \delta} \log(R\mathrm{e}^{\mathrm{i} \theta}) h(x + R\mathrm{e}^{\mathrm{i} \theta}) R \mathrm{e}^{\mathrm{i} \theta} \, \mathrm{d} \theta,%\\
            %&\to& 0 \qquad \text{ as } \varepsilon \to 0 \text{ and then } R\to \infty.
        \end{eqnarray*}
        which tends to zero as $R \to \infty$.
    \end{proof}

\end{document}

%% file: figs/figure1.tex
\definecolor{mycolor1}{rgb}{0.00000,0.44700,0.74100}%
\definecolor{mycolor2}{rgb}{0.85000,0.32500,0.09800}%
\definecolor{mycolor3}{rgb}{0.92900,0.69400,0.12500}%
\definecolor{mycolor4}{rgb}{0.49400,0.18400,0.55600}%
\definecolor{mycolor5}{rgb}{0.46600,0.67400,0.18800}%
\begin{tikzpicture}

\begin{axis}[%
width=5cm,
height=3.6cm,
at={(0cm,0cm)},
scale only axis,
xmin=0,
xmax=40,
xlabel style={font=\color{white!15!black}},
xlabel={$n$ (rank)},
ymode=log,
ymin=1.25e-14,
ymax=10,
yminorticks=true,
ytick={1e-12, 1e-9, 1e-6, 1e-3, 1},
ylabel style={font=\color{white!15!black}},
ylabel={$\|A - A_n \|_2 / \|A\|_2$},
axis background/.style={fill=white},
legend style={at={(1.03,0.5)}, anchor=west, legend cell align=left, align=left, draw=white!15!black}
]
\addplot [color=mycolor1, line width=0.8pt]
  table[row sep=crcr]{%
0	6.27559065214064\\
1	6.21345609122842\\
2	6.03306157675507\\
3	5.75107005852308\\
4	5.39125232024594\\
5	4.98007099220386\\
6	4.542748579824\\
8	3.67006189507278\\
10	2.8841227697305\\
13	1.94995264726971\\
17	1.12692185488243\\
25	0.36654520702383\\
40	0.0440535953497291\\
};
\addlegendentry{Little--Reade \cite{little1984eigenvalues}}

\addplot [color=mycolor1, line width=0.8pt, only marks, mark size=0.8pt, mark=*, mark options={solid, mycolor1}]
  table[row sep=crcr]{%
0	1\\
1	0.31625900530834\\
2	0.191741963997327\\
3	0.14613541204665\\
4	0.120234520316465\\
5	0.102139870367513\\
6	0.0881132201346154\\
7	0.0766073880314548\\
8	0.066860006090983\\
9	0.058445180515523\\
10	0.0510996709702865\\
11	0.0446462723221064\\
12	0.0389569277123568\\
13	0.0339337329048624\\
14	0.0294985487987171\\
15	0.0255869590934844\\
16	0.0221445079765648\\
17	0.0191241789611869\\
18	0.0164845780442651\\
19	0.0141885397074842\\
20	0.0122020106656483\\
21	0.0104931457681274\\
22	0.00903160353177468\\
23	0.00778806743670735\\
24	0.00673404153696101\\
25	0.00584196346221333\\
26	0.00508563300828205\\
27	0.00444087540333847\\
28	0.00388627749473024\\
29	0.00340379940896919\\
30	0.00297910063307283\\
31	0.00260151094928442\\
32	0.00226367543592041\\
33	0.00196096668447141\\
34	0.00169077367939931\\
35	0.00145175854728072\\
36	0.00124314239132788\\
37	0.00106406058444526\\
38	0.000913030932361516\\
39	0.000787606809346363\\
40	0.000684314561852617\\
};
\addlegendentry{Chebyshev polynomial interpolant}

\addplot [color=mycolor2, line width=0.8pt]
  table[row sep=crcr]{%
1	16\\
2	5.78008410230058\\
3	5.78008410230058\\
4	2.08808576435421\\
5	2.08808576435421\\
6	0.754331958173995\\
7	0.754331958173995\\
8	0.272506384956173\\
9	0.272506384956173\\
10	0.0984443639662845\\
11	0.0984443639662845\\
12	0.0355635439451633\\
13	0.0355635439451633\\
14	0.0128475171861817\\
15	0.0128475171861817\\
16	0.00464123311511755\\
17	0.00464123311511755\\
18	0.00167666985898512\\
19	0.00167666985898512\\
20	0.000605705799795404\\
21	0.000605705799795404\\
22	0.000218814404004288\\
23	0.000218814404004288\\
24	7.90478536212226e-05\\
25	7.90478536212226e-05\\
26	2.85564526273133e-05\\
27	2.85564526273133e-05\\
28	1.03161686155771e-05\\
29	1.03161686155771e-05\\
30	3.72677013822176e-06\\
31	3.72677013822176e-06\\
32	1.34631530180401e-06\\
33	1.34631530180401e-06\\
34	4.86363479540076e-07\\
35	4.86363479540076e-07\\
36	1.75701363501824e-07\\
37	1.75701363501824e-07\\
38	6.34730411205893e-08\\
39	6.34730411205893e-08\\
40	2.29299697441118e-08\\
};
\addlegendentry{Beckermann--Townsend \cite{beckermann2019bounds}}

\addplot [color=mycolor3, line width=0.8pt]
  table[row sep=crcr]{%
0	6.27559065214069\\
1	7.39258793440364\\
29	1.00713968666881e-14\\
};
\addlegendentry{Theorem \ref{thm:main}}

\addplot [color=mycolor3, line width=0.8pt, only marks, mark size=1.5pt, mark=o, mark options={solid, mycolor3}]
  table[row sep=crcr]{%
0	1\\
1	0.508820739277728\\
2	0.0735507057274247\\
3	0.0179159547357388\\
4	0.00580851650472821\\
5	0.00163895096186486\\
6	0.000417339000910261\\
7	0.000114658234638617\\
8	3.19968114238311e-05\\
9	8.94664501484438e-06\\
10	2.54245093573951e-06\\
11	7.1148147376949e-07\\
12	1.94750826861245e-07\\
13	6.2336373988455e-08\\
14	2.09442829588703e-08\\
15	5.58569332928419e-09\\
16	1.43687014744987e-09\\
17	4.29077946843053e-10\\
18	1.1958119892331e-10\\
19	3.40470995297498e-11\\
20	1.13054141974073e-11\\
21	3.77613281646295e-12\\
22	1.14058515970803e-12\\
23	2.76093078068248e-13\\
24	7.20582717293295e-14\\
25	3.43235187376533e-14\\
};
\addlegendentry{Zolotarev rational interpolant}

\addplot [color=mycolor4, line width=0.8pt]
  table[row sep=crcr]{%
12	5.31000121910089e-07\\
14	1.66884945420587e-08\\
15	2.5512807847251e-09\\
16	4.10557050188294e-10\\
17	6.37081090583185e-11\\
18	1.01079472037394e-11\\
19	1.4094670360778e-12\\
20	2.03553786588746e-13\\
21	2.80620685138551e-14\\
22	6.81892505540347e-15\\
};
\addlegendentry{Theorem \ref{thm:main} (discrete)}

\addplot [color=mycolor4, line width=0.8pt, only marks, mark=triangle, mark options={solid, mycolor4}]
  table[row sep=crcr]{%
12	1.05712547977934e-08\\
13	1.77260438218948e-09\\
14	3.00770520837656e-10\\
15	4.77070552622078e-11\\
16	7.44868799043584e-12\\
17	1.16003671385335e-12\\
18	1.74000515321254e-13\\
19	2.95408441165873e-14\\
};
\addlegendentry{Discrete Zolotarev rational interpolant}

\addplot [color=mycolor5, line width=0.8pt, mark=asterisk, mark options={solid, mycolor5}]
  table[row sep=crcr]{%
0	1\\
1	0.180397182172517\\
2	0.0547446621921632\\
3	0.0119338932656883\\
4	0.00231401876297875\\
5	0.000419802477097068\\
6	7.18036204758932e-05\\
7	1.16296952396004e-05\\
8	1.7897610836812e-06\\
9	2.62442561048387e-07\\
10	3.67525713634621e-08\\
11	4.92491971721334e-09\\
12	6.32540568078897e-10\\
13	7.79785318868003e-11\\
14	9.2385705722165e-12\\
15	1.05359063283518e-12\\
16	1.15963570502145e-13\\
17	1.26311809281612e-14\\
18	8.9038291207368e-15\\
};
\addlegendentry{Best ($\sigma_{n+1}(A)/\sigma_1(A)$)}

\end{axis}
\end{tikzpicture}%

%% file: figs/figure3.tex
\definecolor{mycolor1}{rgb}{0.00000,0.44700,0.74100}%
\definecolor{mycolor2}{rgb}{0.92900,0.69400,0.12500}%
\definecolor{mycolor3}{rgb}{0.46600,0.67400,0.18800}%
\definecolor{mycolor4}{rgb}{0.49400,0.18400,0.55600}%
\begin{tikzpicture}

\begin{axis}[%
width=4cm,
height=3.6cm,
at={(0cm,0cm)},
scale only axis,
xmin=0,
xmax=25,
xlabel style={font=\color{white!15!black}},
xlabel={$n$ (rank)},
ymode=log,
ymin=1e-14,
ymax=10,
yminorticks=true,
ytick={1e-12, 1e-9, 1e-6, 1e-3, 1},
ylabel style={font=\color{white!15!black}},
ylabel={$\|A - A_n \|_2 / \|A\|_2$},
axis background/.style={fill=white}
]
\addplot [color=mycolor1, line width=0.8pt, forget plot]
  table[row sep=crcr]{%
0	1\\
1	0.942307692307692\\
2	0.798470236115733\\
3	0.626205582405418\\
4	0.467948861042965\\
5	0.340283816654907\\
7	0.17349170726487\\
11	0.0434063607911826\\
25	0.000331792489830831\\
};
\addplot [color=mycolor1, line width=0.8pt, only marks, mark size=0.8pt, mark=*, mark options={solid, mycolor1}, forget plot]
  table[row sep=crcr]{%
0	1\\
1	0.647399120221468\\
2	0.400841399268676\\
3	0.261197532124482\\
4	0.17543872566438\\
5	0.119748939923826\\
6	0.0824407759696659\\
7	0.057046119777216\\
8	0.0396127091605858\\
9	0.0275801981804897\\
10	0.0192416103907063\\
11	0.0134432716450425\\
12	0.00939967448858018\\
13	0.00657318881370848\\
14	0.00459412662354449\\
15	0.00320711269808119\\
16	0.00223492758992389\\
17	0.00155403899514396\\
18	0.00107797812732342\\
19	0.000745995353363495\\
20	0.000515279537573351\\
21	0.000355592587743249\\
22	0.00024554901389713\\
23	0.000170022795123665\\
24	0.000118332252137453\\
25	8.2965297977359e-05\\
};
\addplot [color=mycolor2, line width=0.8pt, forget plot]
  table[row sep=crcr]{%
0	1\\
1	0.641791365424084\\
19	3.18715865834007e-15\\
};
\addplot [color=mycolor2, line width=0.8pt, only marks, mark size=1.5pt, mark=o, mark options={solid, mycolor2}, forget plot]
  table[row sep=crcr]{%
0	1\\
1	0.240386842258584\\
2	0.0523398015514582\\
3	0.00894150438123018\\
4	0.00149120220961286\\
5	0.000243383989006798\\
6	3.85635760589317e-05\\
7	5.97103319526646e-06\\
8	9.32877902905141e-07\\
9	1.51929117590198e-07\\
10	2.55156388646627e-08\\
11	4.18122321019913e-09\\
12	6.43281391673657e-10\\
13	9.9366342041864e-11\\
14	1.63687592414652e-11\\
15	2.68513824584405e-12\\
16	4.23563634536589e-13\\
17	6.69271556078623e-14\\
18	1.20584943551447e-14\\
};
\addplot [color=mycolor3, line width=0.8pt, mark=asterisk, mark options={solid, mycolor3}, forget plot]
  table[row sep=crcr]{%
0	1\\
1	0.235543500657835\\
2	0.0438887532304456\\
3	0.00761787101182415\\
4	0.00126039854343367\\
5	0.000199477667324387\\
6	3.02234331154427e-05\\
7	4.38626593017832e-06\\
8	6.10148730566597e-07\\
9	8.14149510312546e-08\\
10	1.04295310342155e-08\\
11	1.28377146397437e-09\\
12	1.51962039802534e-10\\
13	1.73122315088563e-11\\
14	1.89962260467896e-12\\
15	2.0090034486775e-13\\
16	2.04918408380418e-14\\
17	2.01829675844468e-15\\
};
\end{axis}

\begin{axis}[%
width=4cm,
height=3.6cm,
at={(4.25cm,0cm)},
scale only axis,
xmin=0,
xmax=25,
xlabel style={font=\color{white!15!black}},
xlabel={$n$ (rank)},
ymode=log,
ymin=1e-14,
ymax=10,
yticklabels={\empty},
yminorticks=true,
ytick={1e-12, 1e-9, 1e-6, 1e-3, 1},
axis background/.style={fill=white},
legend style={at={(1.05,0.5)}, anchor=west, legend cell align=left, align=left, draw=white!15!black}
]
\addplot [color=mycolor1, line width=0.8pt]
  table[row sep=crcr]{%
0	1\\
1	0.924528301886794\\
2	0.746347528753495\\
3	0.550407253367266\\
4	0.386034652443316\\
5	0.264063028410418\\
7	0.120030971849763\\
13	0.0108392340250347\\
25	8.77610538599389e-05\\
};
\addlegendentry{Proposition \ref{prop:LR}}

\addplot [color=mycolor1, line width=0.8pt, only marks, mark size=0.8pt, mark=*, mark options={solid, mycolor1}]
  table[row sep=crcr]{%
0	1\\
1	0.381268630013295\\
2	0.18250535044087\\
3	0.104213743796801\\
4	0.0642124170066538\\
5	0.0410221147617906\\
6	0.0267070347263352\\
7	0.0175683423519459\\
8	0.0116183814093815\\
9	0.00769716217831046\\
10	0.00509415605418368\\
11	0.00336047451003925\\
12	0.00220613111046288\\
13	0.00144044511095396\\
14	0.000936177544932377\\
15	0.000607383225871054\\
16	0.000395450001992105\\
17	0.000260171160475098\\
18	0.000174023352120386\\
19	0.00011851416655819\\
20	8.17408940088657e-05\\
21	5.65020936144732e-05\\
22	3.86824323304589e-05\\
23	2.59698273165762e-05\\
24	1.6994158655657e-05\\
25	1.08413004301815e-05\\
};
\addlegendentry{Chebyshev interpolant}

\addplot [color=mycolor2, line width=0.8pt]
  table[row sep=crcr]{%
0	1\\
1	0.702750564104042\\
20	3.13967158568748e-15\\
};
\addlegendentry{Theorem \ref{thm:main}}

\addplot [color=mycolor2, line width=0.8pt, only marks, mark size=1.5pt, mark=o, mark options={solid, mycolor2}]
  table[row sep=crcr]{%
0	1\\
1	0.269549365748745\\
2	0.0548723372830671\\
3	0.011186207514688\\
4	0.0020129134277974\\
5	0.000347961798786948\\
6	6.06828688000521e-05\\
7	1.06767977204534e-05\\
8	1.9046886833679e-06\\
9	3.42397437097339e-07\\
10	5.90294339660345e-08\\
11	9.94040170161492e-09\\
12	1.80981530072354e-09\\
13	3.30777694086859e-10\\
14	5.4947554185512e-11\\
15	9.40285733031044e-12\\
16	1.78621876800294e-12\\
17	3.08983024101558e-13\\
18	5.15092758616877e-14\\
19	1.70597514277102e-14\\
};
\addlegendentry{Zolotarev interpolant}

\addplot [color=mycolor4, line width=0.8pt]
  table[row sep=crcr]{%
13	5.596062731582e-11\\
14	6.15090996343753e-12\\
15	6.16711815152736e-13\\
16	6.03695089321624e-14\\
17	5.53655883817907e-15\\
};
\addlegendentry{Theorem \ref{thm:main} (discrete)}

\addplot [color=mycolor4, line width=0.8pt, only marks, mark=triangle, mark options={solid, mycolor4}]
  table[row sep=crcr]{%
13	3.34527039058277e-11\\
14	3.63636061612976e-12\\
15	3.73884688428412e-13\\
16	3.65104556425706e-14\\
};
\addlegendentry{Discrete Zolotarev interpolant}

\addplot [color=mycolor3, line width=0.8pt, mark=asterisk, mark options={solid, mycolor3}]
  table[row sep=crcr]{%
0	1\\
1	0.192894368432449\\
2	0.042323718021921\\
3	0.0082842705708791\\
4	0.00145513463896555\\
5	0.000238504881411933\\
6	3.75889259728654e-05\\
7	5.65118135592919e-06\\
8	7.92901289071177e-07\\
9	1.04280009065002e-07\\
10	1.31272157852949e-08\\
11	1.59932537677319e-09\\
12	1.86989641966531e-10\\
13	2.07009222398901e-11\\
14	2.18345573144474e-12\\
15	2.21868843759838e-13\\
16	2.14976530851526e-14\\
17	1.97101681050675e-15\\
};
\addlegendentry{Best}

\end{axis}
\end{tikzpicture}%

%% file: figs/figure4.tex
\definecolor{mycolor1}{rgb}{0.00000,0.44700,0.74100}%
\definecolor{mycolor2}{rgb}{0.92900,0.69400,0.12500}%
\definecolor{mycolor3}{rgb}{0.30100,0.74500,0.93300}%
\definecolor{mycolor4}{rgb}{0.46600,0.67400,0.18800}%
\begin{tikzpicture}

\begin{axis}[%
width=5cm,
height=3.6cm,
at={(0cm,0cm)},
scale only axis,
xmin=0,
xmax=40,
xlabel style={font=\color{white!15!black}},
xlabel={$n$ (rank)},
ymode=log,
ymin=1e-14,
ymax=10,
yminorticks=true,
ytick={1e-12, 1e-9, 1e-6, 1e-3, 1},
ylabel style={font=\color{white!15!black}},
ylabel={$\|A - A_n \|_2 / \|A\|_2$},
axis background/.style={fill=white},
legend style={at={(1.03,0.5)}, anchor=west, legend cell align=left, align=left, draw=white!15!black}
]
\addplot [color=mycolor1, line width=0.8pt]
  table[row sep=crcr]{%
1	0.431053594328192\\
2	0.414313648917385\\
3	0.37004084067711\\
4	0.311569713269338\\
5	0.251505088585614\\
6	0.197477470052884\\
7	0.152418660866062\\
9	0.0884234305277412\\
13	0.0287442284968138\\
35	5.65025537022053e-05\\
40	1.370403603805e-05\\
};
\addlegendentry{Proposition \ref{prop:LR}}

\addplot [color=mycolor1, line width=0.8pt, only marks, mark size=0.8pt, mark=*, mark options={solid, mycolor1}]
  table[row sep=crcr]{%
1	0.148097844798918\\
2	0.0333063709206341\\
3	0.0168068228729864\\
4	0.00836336277238304\\
5	0.00440213558940146\\
6	0.00254025921377988\\
7	0.00157202479232702\\
8	0.0010211797456418\\
9	0.000654858394104461\\
10	0.000439375153794305\\
11	0.000290697698084299\\
12	0.000195786285135836\\
13	0.000140497601913563\\
14	9.64727798218682e-05\\
15	6.30558200445427e-05\\
16	4.60036587662307e-05\\
17	3.05863915934713e-05\\
18	2.13461439345386e-05\\
19	1.60548310628823e-05\\
20	1.12700685685549e-05\\
21	7.88390487519172e-06\\
22	5.89150115994294e-06\\
23	4.16706335873364e-06\\
24	3.00613032699178e-06\\
25	2.08349361371777e-06\\
26	1.53644131860581e-06\\
27	1.05747237160028e-06\\
28	7.71564365346415e-07\\
29	5.49122941633122e-07\\
30	4.12635393795286e-07\\
31	2.97630687121614e-07\\
32	2.19069334165389e-07\\
33	1.60554834820275e-07\\
34	1.16979564518552e-07\\
35	8.57935219845555e-08\\
36	6.19433385525611e-08\\
37	4.48578277229011e-08\\
38	3.35445162317113e-08\\
39	2.36020031718172e-08\\
40	1.67653461002422e-08\\
};
\addlegendentry{Chebyshev interpolant}

\addplot [color=mycolor2, line width=0.8pt]
  table[row sep=crcr]{%
0	0.431053594328189\\
1	0.431053594328189\\
2	0.324703786776745\\
3	0.089008444606773\\
5	0.00469549736535038\\
24	3.07917799166547e-15\\
};
\addlegendentry{Equation \eqref{eqn:logCauchybound}}

\addplot [color=mycolor2, line width=0.8pt, only marks, mark size=1.5pt, mark=o, mark options={solid, mycolor2}]
  table[row sep=crcr]{%
1	0.148097844798918\\
2	0.0628327650144977\\
3	0.00486224368804542\\
4	0.0015866537900253\\
5	0.000373017736820345\\
6	6.13414847048066e-05\\
7	8.31350407981782e-06\\
8	1.34172707479323e-06\\
9	3.23581788415849e-07\\
10	5.92828305551223e-08\\
11	1.29310063530242e-08\\
12	2.90047757546412e-09\\
13	6.08536707408595e-10\\
14	1.16949239064214e-10\\
15	2.35326600600609e-11\\
16	5.28970983866001e-12\\
17	1.16801395626744e-12\\
18	2.38405911106825e-13\\
19	6.36202746122464e-14\\
20	2.22770540100871e-14\\
};
\addlegendentry{Zolotarev interpolant}

\addplot [color=mycolor3, line width=0.8pt, only marks, mark=+, mark options={solid, mycolor3}]
  table[row sep=crcr]{%
0	1\\
1	0.634738266589779\\
2	0.175558810874072\\
3	0.053817456552097\\
4	0.0144787926736233\\
5	0.00342842782015129\\
6	0.000858439054586346\\
7	0.000201793887254342\\
8	4.43578413254411e-05\\
9	1.03052139213601e-05\\
10	2.52320805515116e-06\\
11	5.68623267530875e-07\\
12	1.33960436633057e-07\\
13	3.15327498679073e-08\\
14	7.01967214742367e-09\\
15	1.67807074768713e-09\\
16	3.90509065611667e-10\\
17	9.01631185054797e-11\\
18	1.98054982140063e-11\\
19	4.80956660382592e-12\\
20	1.10270099308319e-12\\
21	2.49616225315696e-13\\
22	6.36141927734083e-14\\
23	2.90284391978744e-14\\
};
\addlegendentry{Suboptimal Zolotarev interpolant}

\addplot [color=mycolor4, line width=0.8pt, mark=asterisk, mark options={solid, mycolor4}]
  table[row sep=crcr]{%
0	1\\
1	0.0329368198293267\\
2	0.00377584494340195\\
3	0.000504575385881117\\
4	7.04618318922247e-05\\
5	9.92762526400509e-06\\
6	1.35415572256275e-06\\
7	1.99097449542368e-07\\
8	1.92946789368244e-08\\
9	2.12929266877586e-09\\
10	2.59150222699387e-10\\
11	2.60701390660805e-11\\
12	2.39638945723005e-12\\
13	2.59182129223587e-13\\
14	2.55969844490696e-14\\
15	1.57805501035089e-15\\
};
\addlegendentry{Best}

\end{axis}
\end{tikzpicture}%

%% file: figs/figure5.tex
\definecolor{mycolor1}{rgb}{0.00000,0.44700,0.74100}%
\definecolor{mycolor2}{rgb}{0.92900,0.69400,0.12500}%
\definecolor{mycolor3}{rgb}{0.49400,0.18400,0.55600}%
\definecolor{mycolor4}{rgb}{0.46600,0.67400,0.18800}%
\begin{tikzpicture}

\begin{axis}[%
width=5cm,
height=3.6cm,
at={(0cm,0cm)},
scale only axis,
xmin=0,
xmax=40,
xlabel style={font=\color{white!15!black}},
xlabel={$n$ (rank)},
ymode=log,
ymin=1.25e-14,
ymax=10,
yminorticks=true,
ytick={1e-12, 1e-9, 1e-6, 1e-3, 1},
ylabel style={font=\color{white!15!black}},
ylabel={$\|A - A_n \|_2 / \|A\|_2$},
axis background/.style={fill=white},
legend style={at={(1.03,0.5)}, anchor=west, legend cell align=left, align=left, draw=white!15!black}
]
\addplot [color=mycolor1, line width=0.8pt]
  table[row sep=crcr]{%
1	5.77628067842821\\
2	5.68830178576642\\
3	5.43729278333593\\
4	5.05766212511244\\
5	4.59472908629718\\
6	4.09343002412052\\
7	3.59061361798238\\
9	2.67312627934921\\
11	1.93674799654868\\
14	1.16542599191455\\
20	0.409935657973778\\
40	0.0122338566238708\\
};
\addlegendentry{Proposition \ref{prop:LR}}

\addplot [color=mycolor1, line width=0.8pt, only marks, mark size=0.8pt, mark=*, mark options={solid, mycolor1}]
  table[row sep=crcr]{%
1	0.854603372658628\\
2	0.466188755372138\\
3	0.423071042186668\\
4	0.321667599901918\\
5	0.246113792087575\\
6	0.190727539607033\\
7	0.149446159694157\\
8	0.118118930264615\\
9	0.0939776135460507\\
10	0.0751408844868653\\
11	0.0602963769485919\\
12	0.0485053778090221\\
13	0.0390815877765757\\
14	0.0315138862125629\\
15	0.0254156230309923\\
16	0.0204903337541916\\
17	0.0165079459060234\\
18	0.0132879022438494\\
19	0.0106869779161332\\
20	0.00859036360164894\\
21	0.00690507048941701\\
22	0.00555501969299638\\
23	0.00447738163212587\\
24	0.00361986870838341\\
25	0.00293877462672585\\
26	0.00239760125221068\\
27	0.00196612327403892\\
28	0.00161972831591476\\
29	0.00133886645696678\\
30	0.00110847673540238\\
31	0.000917329056088225\\
32	0.000757296633187465\\
33	0.000622621319255953\\
34	0.000509240179437887\\
35	0.000414219440653405\\
36	0.000335312260050025\\
37	0.000270633967153025\\
38	0.00021843821382366\\
39	0.000176979944001512\\
40	0.000144461940188072\\
};
\addlegendentry{Chebyshev interpolant}

\addplot [color=mycolor2, line width=0.8pt]
  table[row sep=crcr]{%
0	5.77628067842821\\
1	5.77628067842821\\
2	4.84577564955352\\
3	1.70463405343429\\
4	0.478144578512502\\
12	1.55321142321191e-05\\
29	4.48690477173111e-15\\
};
\addlegendentry{Equation \eqref{eqn:Hankelbound}}

\addplot [color=mycolor2, line width=0.8pt, only marks, mark size=1.5pt, mark=o, mark options={solid, mycolor2}]
  table[row sep=crcr]{%
1	0.854603372658621\\
2	0.216563400022647\\
3	0.0704583818935889\\
4	0.0116492675819745\\
5	0.00274616455507912\\
6	0.000671973168213644\\
7	0.000165969750686502\\
8	4.13580973582632e-05\\
9	1.05753650491098e-05\\
10	2.5859186523307e-06\\
11	6.50815531811036e-07\\
12	1.67738020700163e-07\\
13	4.42276822154142e-08\\
14	1.11380439828822e-08\\
15	2.92757379650274e-09\\
16	7.8455151445371e-10\\
17	2.09799852363606e-10\\
18	5.3532417602816e-11\\
19	1.44742801120093e-11\\
20	3.83411160343891e-12\\
21	1.01630377151172e-12\\
22	2.71146229667815e-13\\
23	7.68036796497286e-14\\
};
\addlegendentry{Zolotarev interpolant}

\addplot [color=mycolor3, line width=0.8pt, only marks, mark=triangle, mark options={solid, mycolor3}]
  table[row sep=crcr]{%
12	1.03409908688126e-08\\
13	1.92719371588328e-09\\
14	3.48284254185416e-10\\
15	6.25086019627658e-11\\
16	1.09677681981764e-11\\
17	1.92655096946667e-12\\
18	3.2384703561612e-13\\
19	5.53266917545115e-14\\
};
\addlegendentry{Semidiscrete Zolotarev interpolant}

\addplot [color=mycolor4, line width=0.8pt, mark=asterisk, mark options={solid, mycolor4}]
  table[row sep=crcr]{%
0	1\\
1	0.0963888746940918\\
2	0.0106192660200535\\
3	0.00114646586848608\\
4	0.000117594768343236\\
5	1.13720222487776e-05\\
6	1.03619232654976e-06\\
7	8.90876216817837e-08\\
8	7.24140575274888e-09\\
9	5.57614164777593e-10\\
10	4.07549113848605e-11\\
11	2.83220104635031e-12\\
12	1.87438947825782e-13\\
13	1.18269355367332e-14\\
};
\addlegendentry{Best}

\end{axis}
\end{tikzpicture}%